\documentclass{amsart}
\usepackage{comment}
\usepackage{mathtools}
\usepackage{placeins}
\usepackage{hyperref}

\mathtoolsset{showonlyrefs=true}

\newtheorem{theorem}{Theorem}

\newtheorem{lemma}[theorem]{Lemma}

\title{Revisiting de Moivre--Laplace}
\author{Raphaël Cerf
}
\address{\kern-8pt Universit\'e Paris-Saclay, CNRS, Laboratoire de math\'ematiques d'Orsay, 91405, Orsay.}
\makeatletter
\renewcommand{\footnoterule}{
  \kern 9pt
  \hrule width 0.4\columnwidth
  \kern 2pt
}
\makeatother
\begin{document}
\begin{abstract}
We revisit the proof of the de Moivre--Laplace theorem, 
which is the ancestor of the central limit theorem for 
the binomial distribution.
Our goal is to provide a proof that can be reasonably presented to 
undergraduate students within a basic course of probability theory.
We follow the strategies presented in two classical references,
the books of Breiman and Feller, but we replace the arguments 
involving series expansions of the logarithm or the exponential 
by the basic inequality 
$\exp(t)\geq
1+t$.
This way we
avoid completely the use of uniform convergence and power series.
We also avoid using Stirling's formula, instead 
we use the exact formula for the Wallis integral.
As a by product of the proof, we also obtain 
a non-asymptotic inequality linking 
the binomial and the Gaussian distributions.
\end{abstract}
\maketitle
\begin{comment}
\thispagestyle{plain}
\footnotetext{
\noindent 
%\footnote{This work was one of the goals of the ANR Project PPPP, number ANR-16-CE40-0016.}
\noindent
Universit\'e Paris-Saclay, CNRS, Laboratoire de math\'ematiques d'Orsay, 91405, Orsay.}

%\begingroup
%\renewcommand\thefootnote{}     % vide le numéro
%%\footnotetext{Keywords: central limit theorem, binomial distribution}
%%\footnotetext{Email: Raphael.Cerf@gmail.com}
%\endgroup
\end{comment}

\def\zk{\{\,1,\dots,k\,\}}
\def\zu{\{0,1\}}
\def\ud{\{\,1,\dots,d\,\}}
\def\unn{\{\,1,\dots,n\,\}}
\def\unns{\{\,1,\dots,n/2\,\}}
\def\zun{\{\,0,\dots,n\,\}}
\def\zdn{\{\,0,\dots,2n\,\}}
\def\zdnp{\{\,0,\dots,2n+1\,\}}
\def\zmn{\{\,-n,\dots,n\,\}}
\def\omegab{\overline{\omega}}
\def\rcurs{\text{\begin{cursive}rl\end{cursive}}}
 \def \Z {{\mathbb Z}}
 \def \E {{\mathbb E}}
 \def \Zd {{\mathbb Z}^d}
 \def \Sd {S^{d-1}}
\def \wQ {\widetilde{Q}}
 \def \R {{\mathbb R}}
 \def \Rd {{\mathbb R}^d}
 \def \C {{\mathbb C}}
 \def \cD {{\mathcal D}}
 \def \tX {\smash{\widetilde{X}}}
 \def \cE {{\mathcal E}}
 \def \cV {{\mathcal V}}
 \def \cH {{\mathcal H}}
 \def \cR {{\mathcal R}}
 \def \cF {{\mathcal F}}
 \def \cC {{\mathcal C}}
 \def \cN {{\mathcal N}}
 \def \cT {{\mathcal T}}
 \def \tT {{\mathcal T}}
 \def \tpi {{\widetilde\pi}}
 \def \cP {{\mathcal P}}
 \def \cPl {{\mathcal P}_\ell}
 \def \cS {{\mathcal S}}
 \def \tS {\smash{\widetilde S}}
 \def \bC {{\overline C}}
 \def \bE {{\overline E}}
 \def \bA {{\overline A}}
 \def \bB {{\overline B}}
 \def \ta {\text{two--arms}}
 \def \N {{\mathbb N}}
 \def \P {{\mathbb P}}
 \def \E {{\mathbb E}}

\newcommand{\set}[2]{\big\{\, #1 :  #2\, \big\}}
\newcommand{\setscolon}[2]{\big\{ #1\, ; \, #2 \big\}}
\newcommand{\Set}[2]{\Big\{\, #1\, : \, #2 \,\Big\}}

\newcommand{\bd}{\partial\, } 
\newcommand{\din}{\partial^{\, in}}
\newcommand{\dout}{\partial^{\, out}}
\newcommand{\dini}{\partial^{\, in}_\infty}
\newcommand{\douti}{\partial^{\, out}_*}
\newcommand{\dine}{\partial^{\, in}_{ext}}
\newcommand{\doutie}{\partial^{\, out,ext}_{\infty}}
\newcommand{\doute}{\partial^{\, out,ext}}
\newcommand{\dedge}{\partial^{\, edge}}
\newcommand{\dexte}{\partial^{ext,\, edge}}
\newcommand{\dstar}{\partial^{\, *}}
\newcommand{\dcirc}{\partial^{\, o}}
\newcommand{\diam}{\text{diameter}\,}
\newcommand{\shell}{\text{Shell}\,}

\newcommand{\ii}{{\underline{i}}} 
\newcommand{\uu}{{\underline{u}}} 
\newcommand{\jj}{{\underline{j}}} 
\newcommand{\kk}{{\underline{k}}} 
\renewcommand{\ll}{{\underline{l}}} 
\newcommand{\mm}{{\underline{m}}} 
\newcommand{\nn}{{\underline{n}}} 
\newcommand{\xx}{{\underline{x}}} 
\newcommand{\zz}{{\underline{z}}} 
\newcommand{\yy}{{\underline{y}}} 
\newcommand{\oo}{{\underline{0}}} 
\newcommand{\rr}{{\underline{r}}} 
\newcommand{\iinf}{{\underline{$\infty$}}} 

\renewcommand{\AA}{\underline{A}}
\newcommand{\EE}{\underline{E}}
\newcommand{\OO}{\underline{O}}
\newcommand{\BB}{\underline{B}}
\newcommand{\DD}{\underline{D}}
\newcommand{\CC}{\underline{C}}
\newcommand{\cCC}{\underline{\calC}}
\newcommand{\FF}{\underline{F}}
\newcommand{\HH}{\underline{H}}
\newcommand{\Su}{\underline{S}} 
\newcommand{\MM}{\underline{M}}  
\newcommand{\RR}{\underline{R}}

\newcommand{\tQ}{\smash{\widetilde{Q}}}
\newcommand{\La}{\Lambda}
\newcommand{\Lak}{\La(k)} 
\newcommand{\Lan}{\La(n)}
\newcommand{\Lam}{\La(m)}
\newcommand{\LaN}{\La(N)}
\newcommand{\LaM}{\La(M)}
\newcommand{\Lann}{{\La(n-\sqrt n)}} 
\newcommand{\Lanf}{\La(n,\phi)}
\newcommand{\LLanf}{\underline{\La}(n,\phi)}
\newcommand{\LLan}{\underline{\La}(n)}
\newcommand{\LLam}{\underline{\La}(m)}
\newcommand{\LLa}{\underline{\La}}
\newcommand{\M}{{\mathbb{M}}} 
\newcommand{\T}{{\mathbb{T}}} 
\newcommand{\Td}{{\mathbb{T}^d}} 
\newcommand{\Eb}{{\mathbb{E}}} 
\renewcommand{\L}{{\mathbb{L}}} 
\renewcommand{\H}{{\mathbb{H}}} 
\newcommand{\Zdn}{{\mathbb{Z}}^d_n} 
\newcommand{\ZZd}{{\underline{\mathbb{Z}}^d}} 
\newcommand{\oomega}{{\underline{\omega}}} 
\newcommand{\ze}{\underline{0}}
\newcommand{\Hd}{{\mathbb{H}}^d}
\newcommand{\Dd}{{\mathbb{D}}^d} 
\newcommand{\Nd}{{\mathbb{N}}^d}
\newcommand{\Ld}{{\mathbb{L}}^d}
\newcommand{\Ldi}{\smash{{\mathbb{L}}^{d,\infty}}}
\newcommand{\Edi}{{\mathbb{E}}^{d,\infty}}
\newcommand{\Ldstar}{\cal{L}^{d,*}} 
\newcommand{\Ed}{{\mathbb{E}}^d} 
\newcommand{\Edstar}{\cal{E}^{d,*}} 
\newcommand{\Eds}{{\mathbb{E}}^{d,*}}

\newcommand{\dist}{{\rm dist}} 
\newcommand{\distT}{\, {\rm dist}_{\calF}} 
\newcommand{\distH}{\, {\rm dist}_{\calH}} 
\newcommand{\DistP}{\, {\rm dist}} 
\newcommand{\dinf}{ {\rm d}_{\infty} }
\newcommand{\di}{{\rm d}_\infty\, }
\newcommand{\done}{\, {\rm d}_{1}\, }
\newcommand{\dtwo}{\, {\rm d}_{2}\, }

%\vfill\eject
%\section{Main results}
%\vspace{-0.75cm}
%\vspace{-0.25cm}
\section{Introduction.} \label{MRF}
%I am in charge of an introductory course on Probability at my University.
I am in charge of an introductory course on Probability at University Paris--Saclay,
for second-year undergraduate students.
The course focuses solely on discrete probability, 
%but my secret dream is to be able to show the students the central limit theorem, 
but my dream is to show the students the central limit theorem, 
or rather a simplified version of it. 
This would help them understand the need to develop the continuous side of 
probability theory, and would serve as an excellent bridge before entering 
the more advanced probability course involving distribution with densities,
and possibly measure theory. 
\vspace{0.25cm}
%\centerline{
%%$B(n,p)$, $p=\frac12$,
%$B(n,\frac12)$, $n=1,2,3,4,5,6,7,8,20,40,50,100$
%}
%%$p=1/2$, $n=1,2,3,4,5,6,7,8,20,40,50,100$ :
%Nous observons que, 
%Lorsque $n$ grandit, nous voyons que 
%Lorsque $n\to\infty$,
%Lorsque $n$ tend vers $+\infty$,
\begin{figure}[ht]
\centering
%\newcolumntype{?}{!{\vrule width 1pt}}
%\vskip-20pt
%\setlength{\arrayrulewidth}{4pt}
\includegraphics[width=0.24\textwidth]{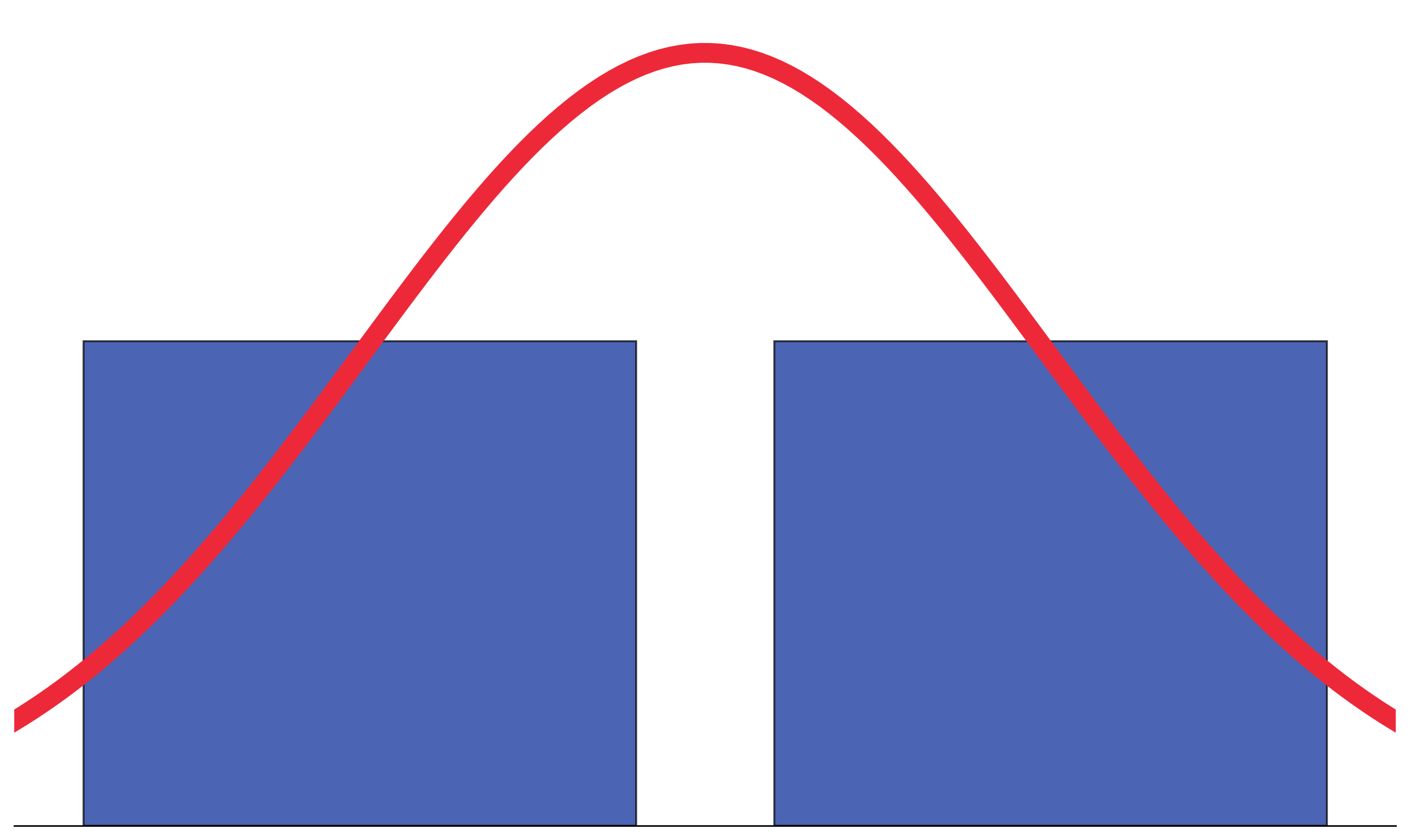}
\kern-3pt
\includegraphics[width=0.24\textwidth]{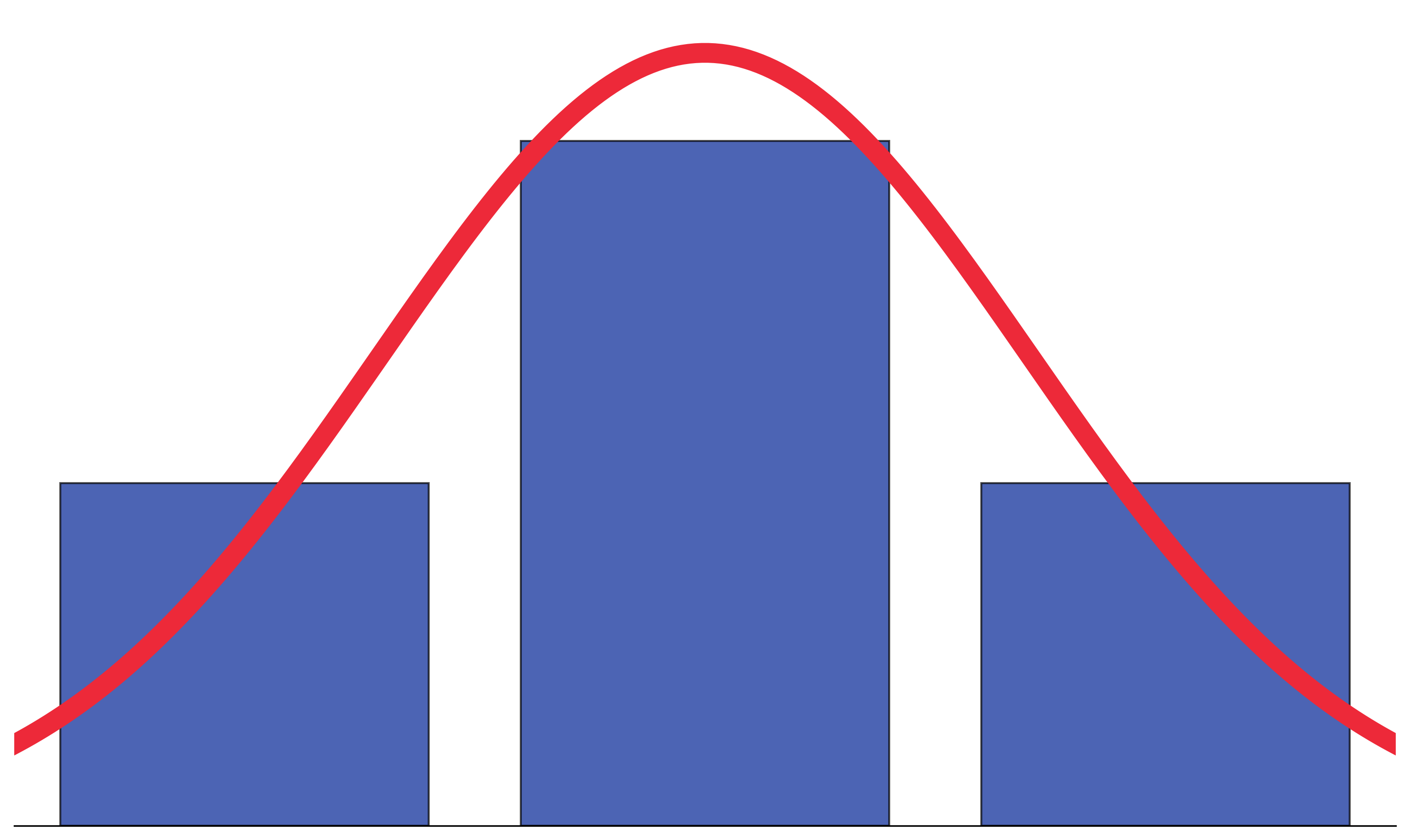}
\kern-3pt
\includegraphics[width=0.24\textwidth]{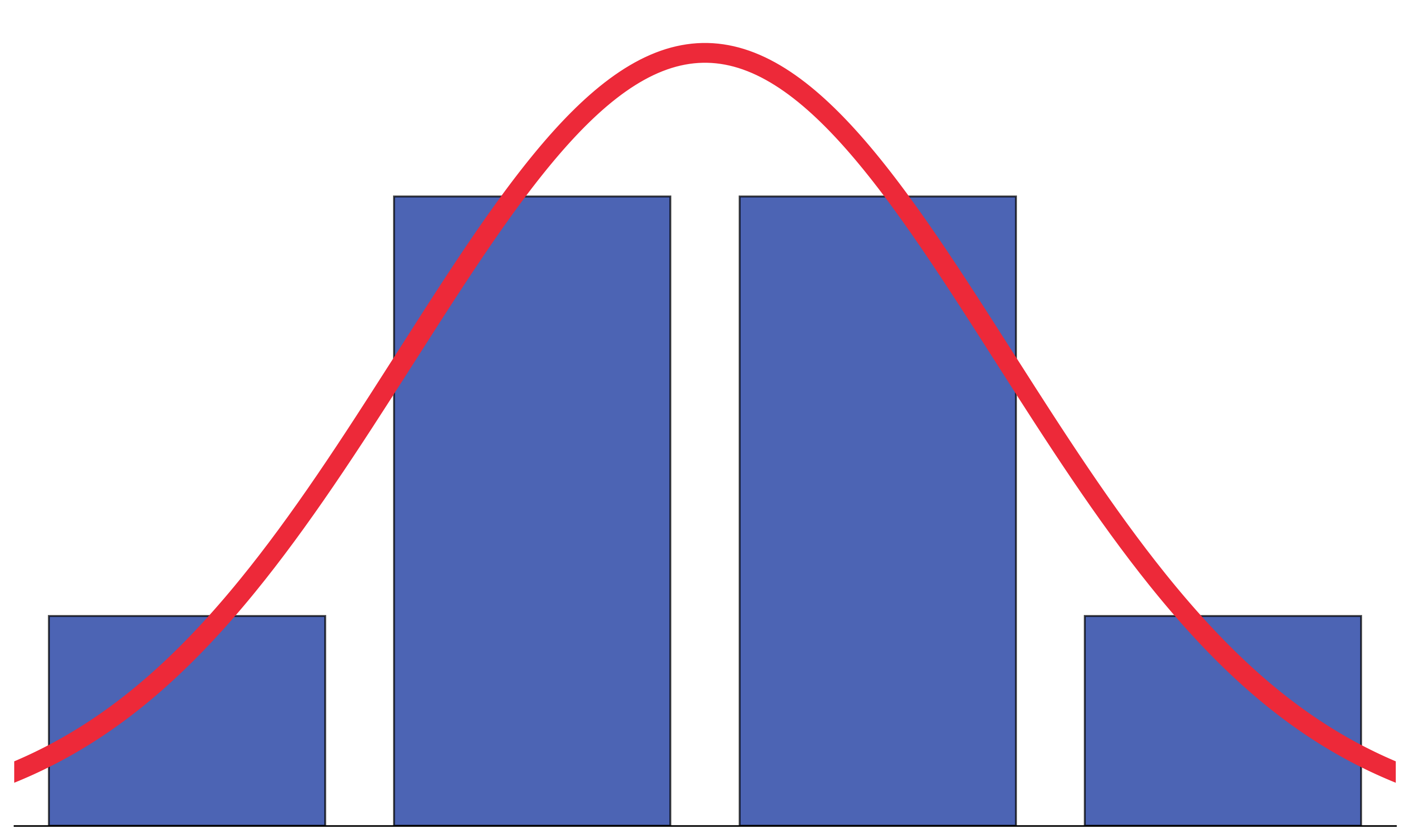}
\kern-3pt
\includegraphics[width=0.24\textwidth]{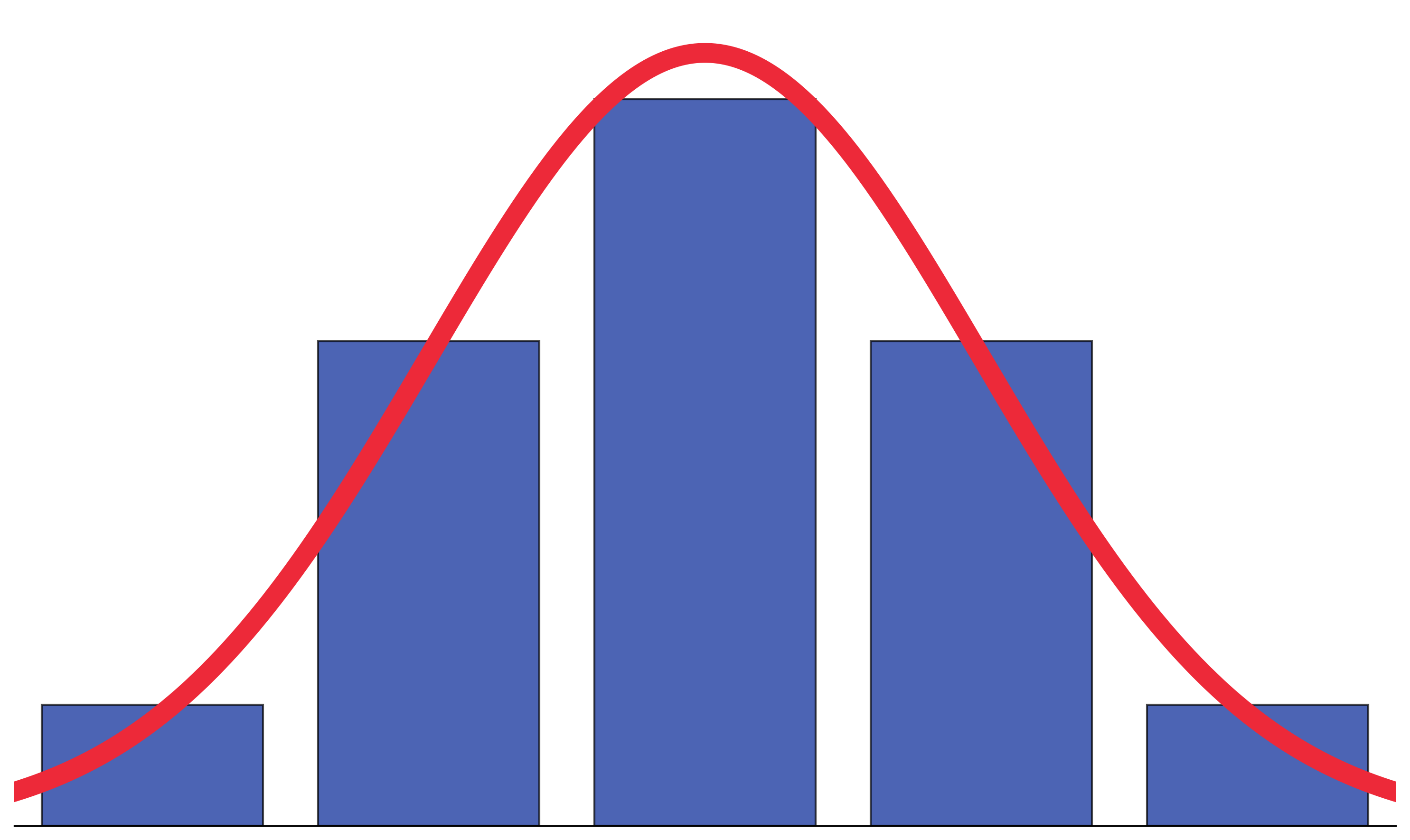}
\bigskip
\smallskip

\includegraphics[width=0.24\textwidth]{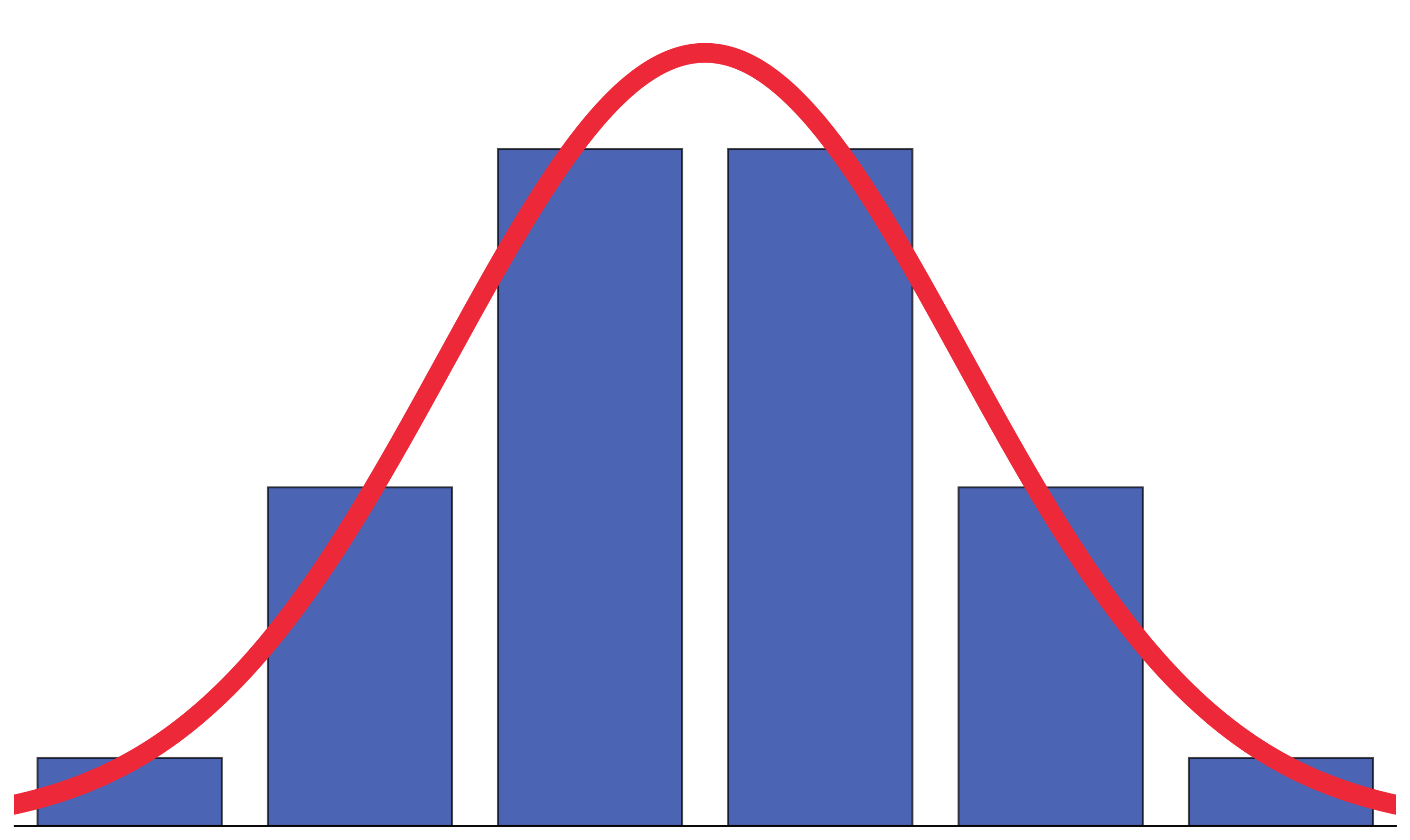}
\kern-3pt
\includegraphics[width=0.24\textwidth]{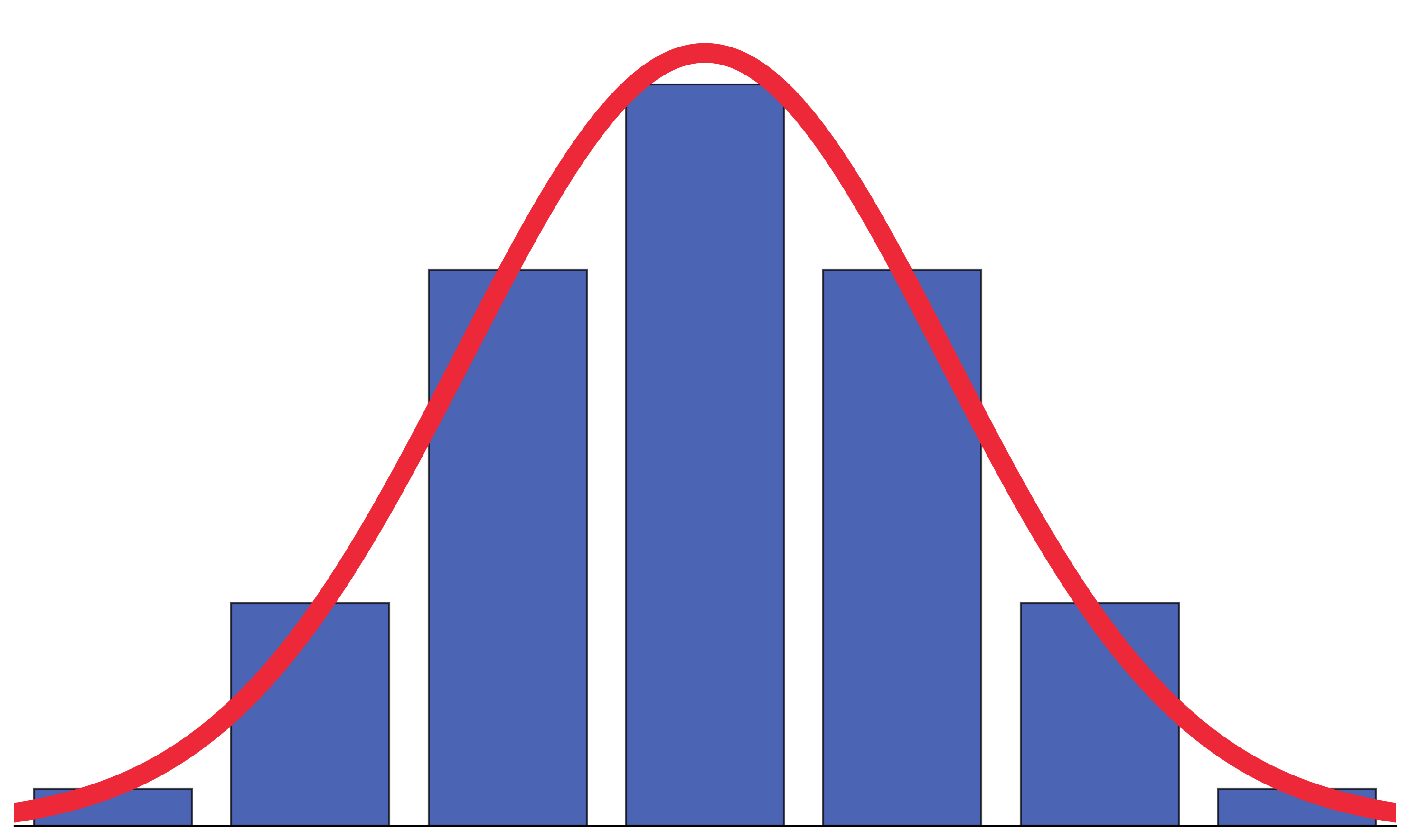}
\kern-3pt
\includegraphics[width=0.24\textwidth]{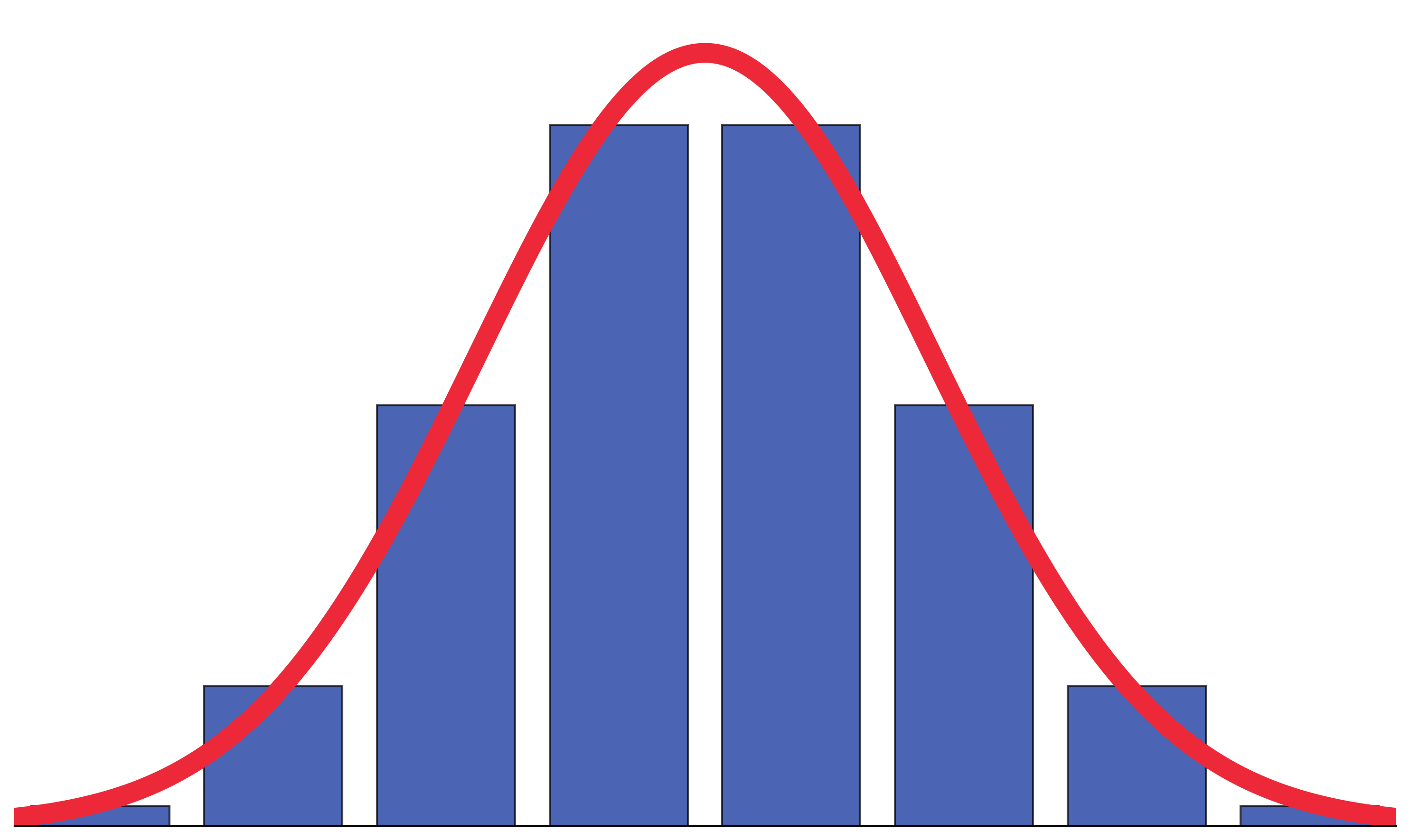}
\kern-3pt
\includegraphics[width=0.24\textwidth]{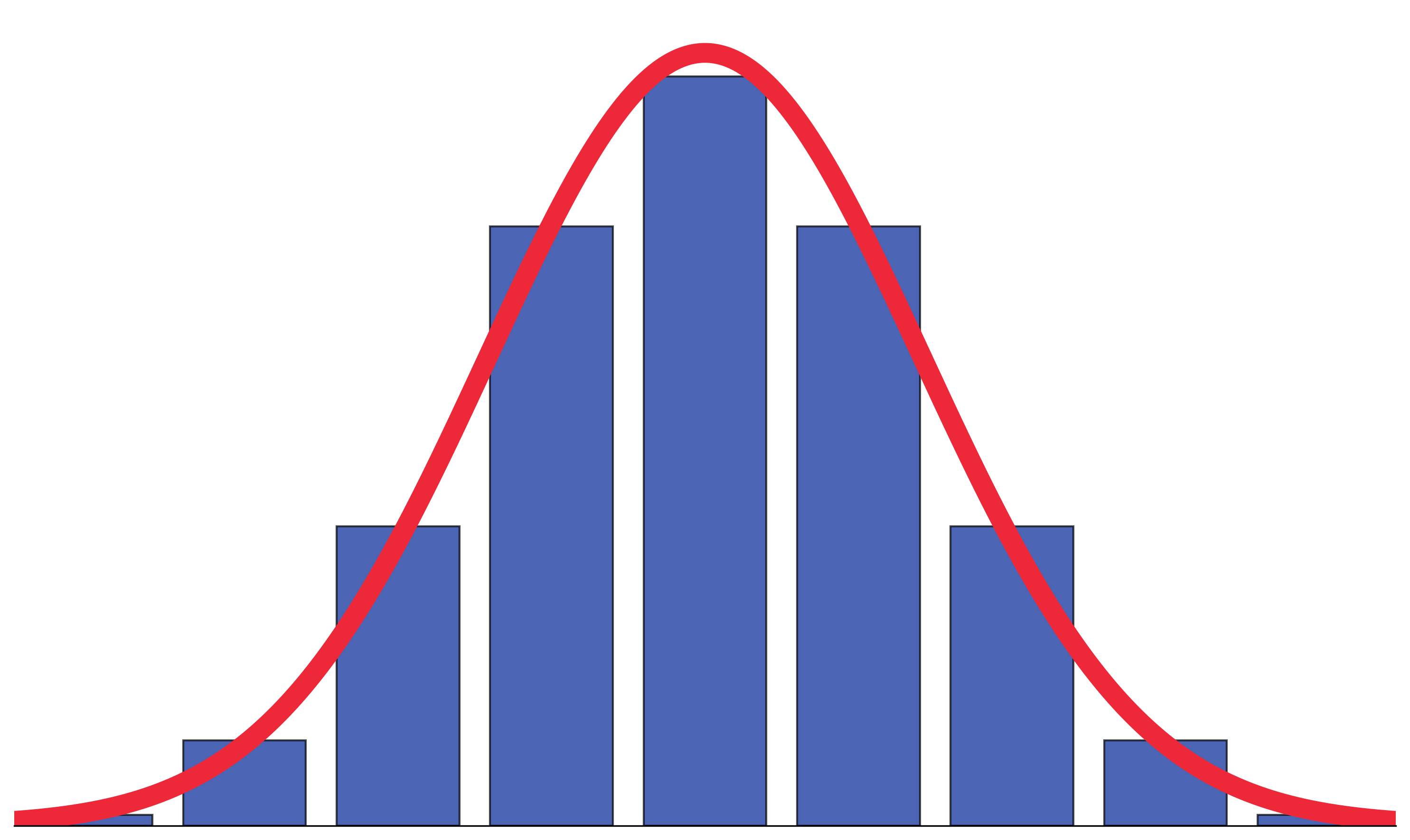}
\bigskip
\smallskip

\includegraphics[width=0.24\textwidth]{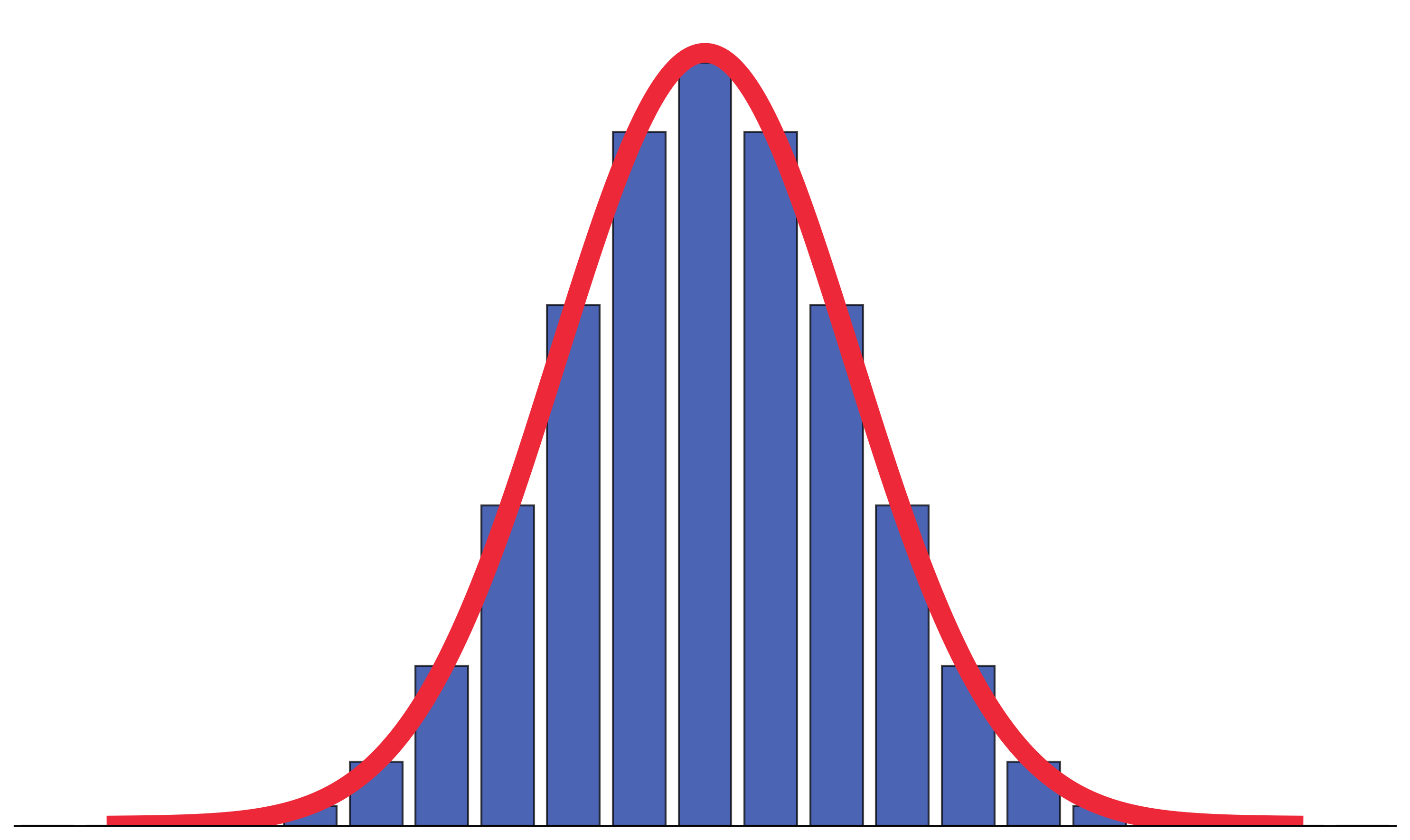}
\kern-3pt
\includegraphics[width=0.24\textwidth]{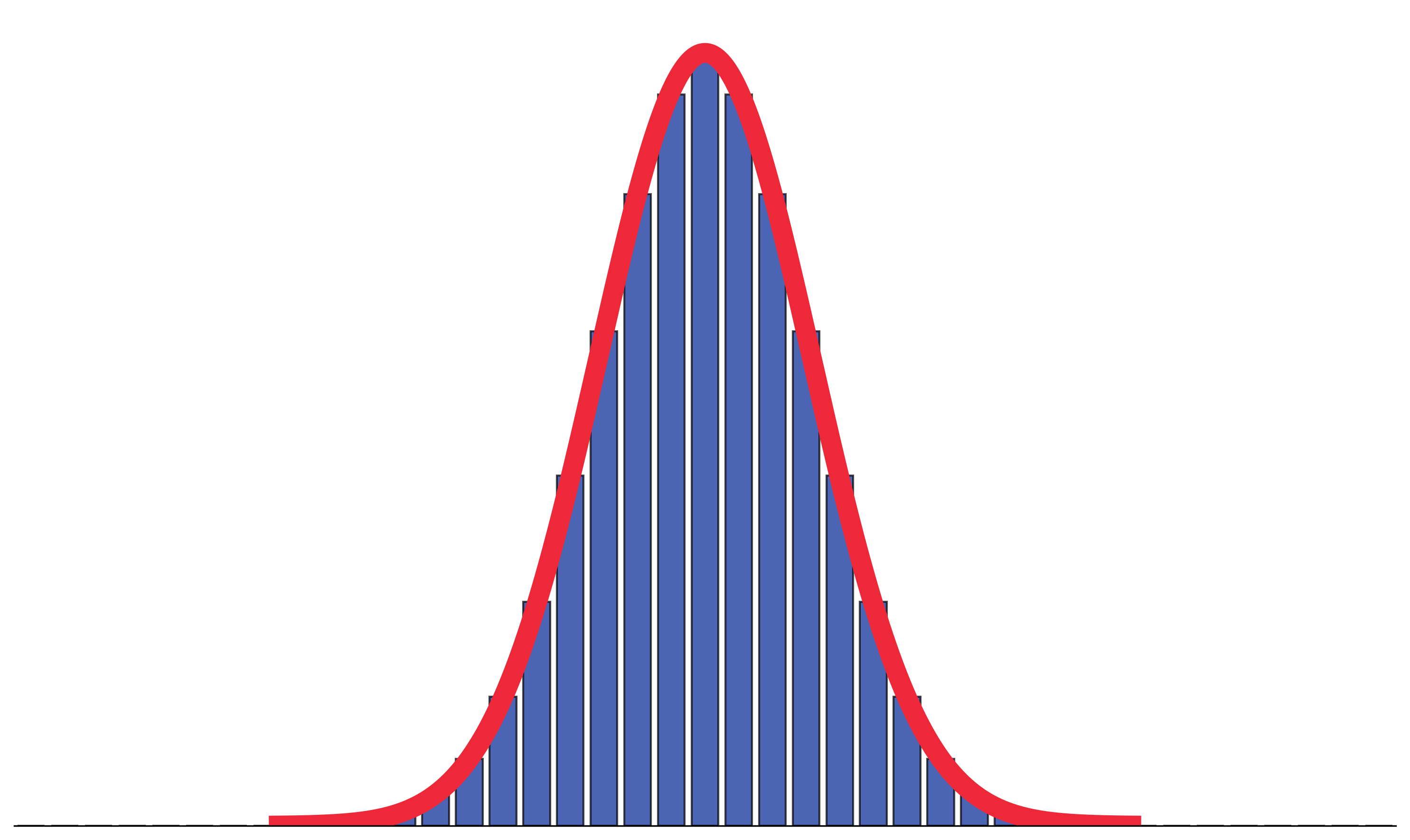}
\kern-3pt
\includegraphics[width=0.24\textwidth]{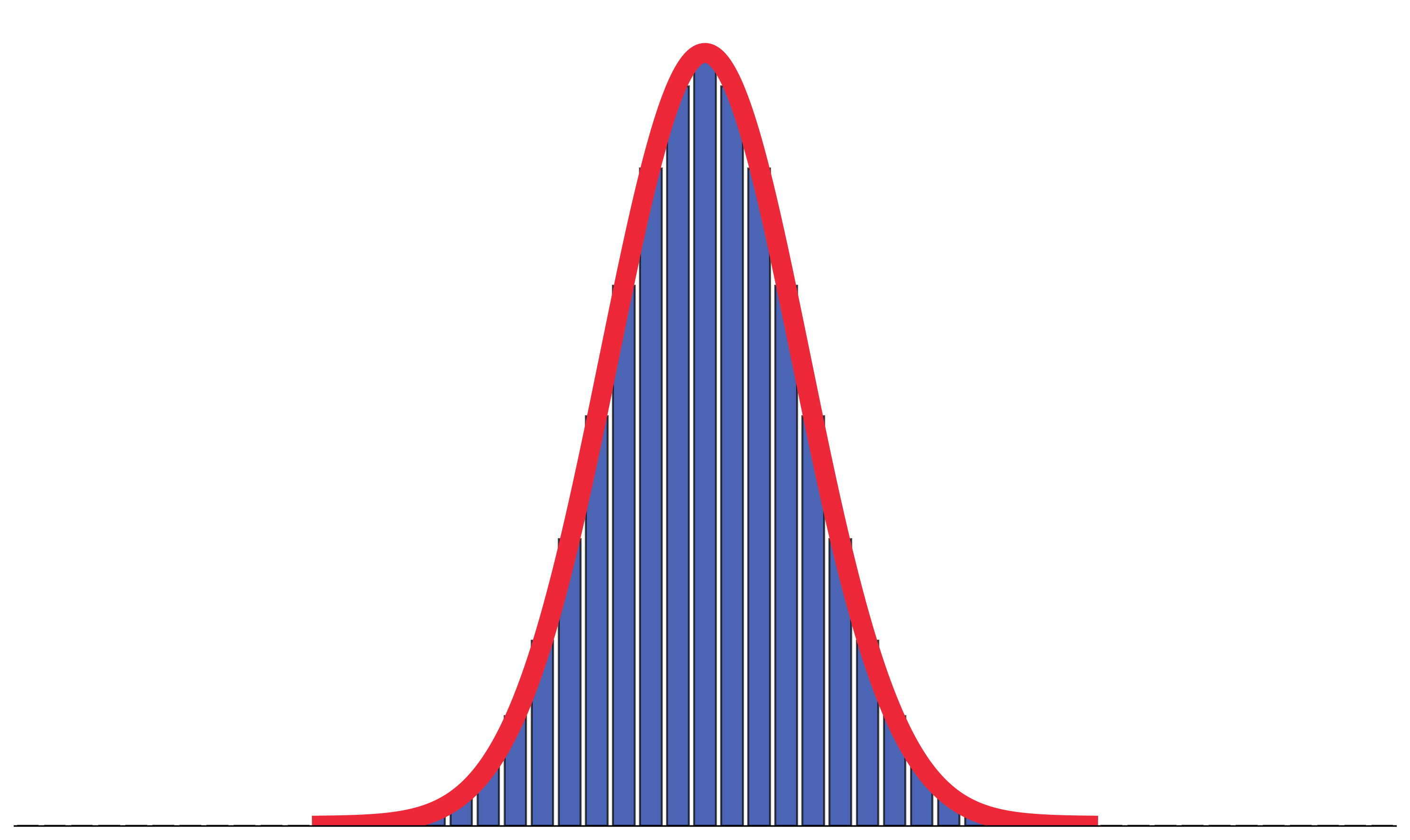}
\kern-3pt
\includegraphics[width=0.24\textwidth]{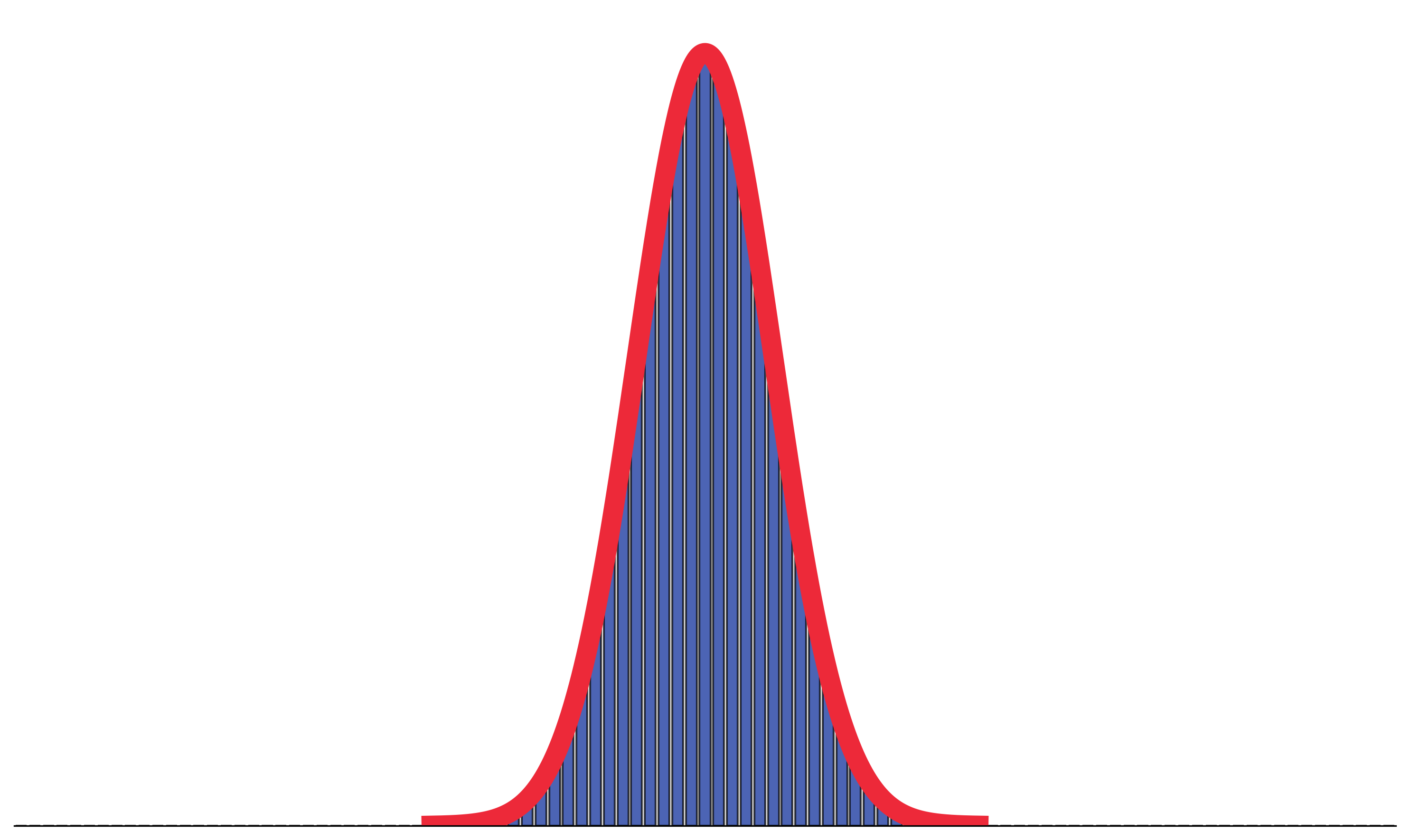}
\setlength\tabcolsep{0 pt}
\renewcommand{\arraystretch}{0}
\begin{comment}{ccc}
\includegraphics[width=0.33\textwidth]{table/binomial_n40_p0.500.png}
	& 
\includegraphics[width=0.33\textwidth]{table/binomial_n50_p0.500.png}
	& 
\includegraphics[width=0.33\textwidth]{table/binomial_n100_p0.500.png}
\end{comment}
%\vskip-10pt
%\caption{The set $A$ is $\text{circle}(0,256)$.
	%The set $D$ is $\Lambda(1024)$.
\smallskip
\caption{
  The laws
$B(n,\frac12)$, 
$\cN(\frac{n}{2},\frac{n}{4})$
%and the corresponding Gaussian distributions, 
for $n=1,2,3,4,5,6,7,8,20,40,50,100$
}
\label{histobig}
	%\vbox{ \hbox{the set $C(A,\Lambda)$ for
%\vskip-50pt
\end{figure}
\FloatBarrier
\noindent
%\newpage
Of course, it is out of the question to aim for the general central limit theorem.
A sensible goal would be to prove the de Moivre--Laplace theorem, which is the 
ancestor of the central limit theorem for the binomial distribution.
It was obtained in the symmetric case $p=\frac12$ by Abraham de Moivre in 1733
and it appeared 
in the second edition of The Doctrine of Chances, 
published in 1738 \cite{M}. 
It is illustrated in figure~\ref{histobig}, which presents
the histograms of the first binomial distributions.
Laplace 
%\cite{L}
proved the result in 1812 
for any  value of $p$ in $]0,1[$
(see chapter III of \cite{L}).
%by Abraham de Moivre
Here is the theorem, in all its splendor, presented in a modern form.
\begin{theorem}[de Moivre--Laplace]
  \label{dml}
  Let $p\in]0,1[$.
  For each $n\geq 1$, let 
$X_1,\dots,X_n$ be $n$ independent Bernoulli random variables with parameter $p$, i.e.,
$$\forall i\in\unn\qquad 
P(X_i=0)\,=\,1-p\,,\qquad
P(X_i=1)\,=\,p\,.$$
%The random variable $X_1+\cdots+X_n$ con
For any real numbers $a<b$, the following limit holds:
%we have the limit
\begin{equation}
  \lim_{n\to\infty}\,P\Big(
    %a<\frac{1}{n}(X_1+\cdots+X_n)<b
    a<\frac{X_1+\cdots+X_n-np}{\sqrt{np(1-p)}}<b
  \Big) \,=\,\frac{1}{\sqrt{2\pi}}
  \int_a^b
   e^{- \tfrac{s^2}{2} } ds\,.
  %\exp\big(-\frac{x^2}{2}\big)\,dx\,.
  %\int_a^b\exp\big(-\frac{x^2}{2}\big)\,dx\,.
\end{equation}
%We provide a detailed proof of Hoeffding's inequality in this case.
%The proof is short (one page) and elementary, but it is quite tricky. 
%\vspace{3.5cm}
\end{theorem}
\vspace{-0.7cm}
%\vspace{-1.5cm}
\begin{figure}[ht]
%\vspace{-0.5cm}
\centering
\includegraphics[width=0.8\textwidth]{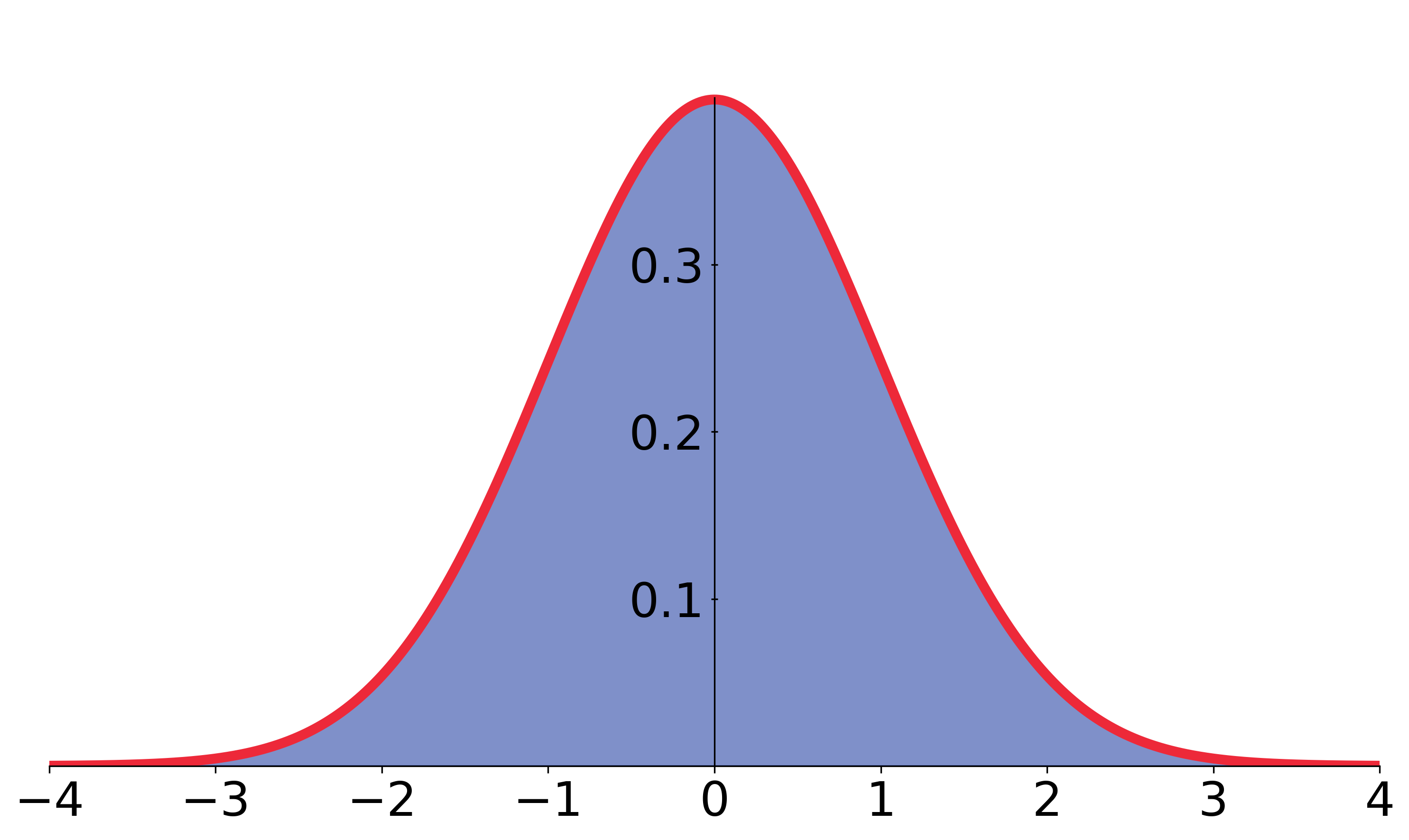}
%\newcolumntype{?}{!{\vrule width 1pt}}
\vspace{-0.2cm}
\caption{
  The Gaussian distribution $\cN(0,1)$
}
\label{gistobig}
	%\vbox{ \hbox{the set $C(A,\Lambda)$ for
%\vskip-50pt
\end{figure}So I searched the literature for an elementary exposition of 
the de Moivre--Laplace theorem. 
%There is a Wikipedia page for the de Moivre--Laplace theorem. It presents the sketches of 
The Wikipedia page for the de Moivre--Laplace theorem presents the sketches of 
two possible arguments. 
The first one is based on a differential equation, which can 
be made rigorous 
(see \cite{JR}), but is not very intuitive.
%The second one follows more the classical approach and relies on Stirling's formula.
The second one follows more the historical approach and relies on Stirling's formula.
Amazingly enough, this theorem, 
which is almost 300 years old,
%dates back to the eighteenth century,
%which dates back to the eighteenth century,
%and appeared in the second edition of The Doctrine of Chances by Abraham de Moivre, published in 1738 \cite{M}, 
%appeared for the first time in the historical treaty of Abraham de Moivre 
is still the subject of recent papers \cite{AG,HM} (these are technical refinements and they are 
not suited for preparing an undergraduate lecture).

The best sources I found were in fact contained in two classical books
of probability. 
Leo Breiman \cite{B} presents a three pages proof in the section 
entitled "The bell-shaped curve enters".
William Feller,
in his legendary book \cite{F}, devotes a full chapter to the normal approximation
of the binomial distribution.
%, in which there are two subsections devoted to the De Moivre--Laplace theorem.
These are the most natural and obvious places to look, and indeed these two 
references provide a full and elementary derivation. Both references present 
first the specific case of the symmetric binomial distribution, 
where the formulas are simpler and the computations easier.

So I wondered: would it be a pleasant and profitable experience for the students
to be exposed to these proofs ?
Profitable, certainly, they would learn a lot from it. Pleasant, I am not so sure,
because the proofs involve technical computations, which seem unavoidable. 
%Let us discuss the passages that might loo
Let us discuss the passages that are likely to leave 
many of our students behind, however willing they may be.

\noindent
$\bullet$ 
The most off-putting part of Breiman's proof is the use of a uniform expansion of the logarithm.
Throughout the calculation, a residual term appears, denoted 
successively by 
$\varepsilon_{j,n}$,
$\varepsilon'_{j,n}$,
$\Delta_{j,n}$, where the index $j$ varies in a set of integers $R_n$. 
The point is that these error terms vanish uniformly over $j$ in $R_n$ 
as $n$ goes to $\infty$. 
The concern is that our students, 
at the beginning of their second year, have not yet mastered uniform convergence.

\noindent
$\bullet$ 
Let us look at Feller's book, 
%3rd edition, 
chapter VII, section 2: "Orientation: Symmetric distribution".
The proof seems less technical than the one of Breiman, with apparently no nasty computations.
Yet, upon further examination, a crucial step is the approximation presented on formula~(2.3) page 180 of \cite{F},
%there,
\begin{equation}
  \label{appro}
1+\frac{j}{\nu}\,=\,e^{\frac{j}{\nu}+\cdots}\,,
\end{equation}
where it is said that ``the dots indicate terms which add to less than 
${j}/{\nu^2}$."
%\og texte cité \fg{}
The rigorous justification comes at the end of the section.
%it involves the series expansion of the logarithm,
It relies on the series expansion of the logarithm, 
%and the comparison of these series 
which is compared with an adequate geometric 
series, in order to
control 
the error in~\eqref{appro}, uniformly over $j$ much smaller than $\nu$.
The problem is that our students 
have not yet learned the series expansions of classical functions.
At this point, it seems that the experience would be very unpleasant, 
not only 
for the students but also for the teacher,
who would find himself explaining a delicate new tool in the middle of an 
already complicated proof.
Hence, our question remains unresolved :
Is it reasonable and possible to present a proof of the de Moivre--Laplace theorem
in the symmetric case to students who do not know series expansion nor uniform convergence ?

Our purpose here is to show that the answer is in fact positive. 
We present a 
%full proof, 
%which is 
variant of the proofs appearing in the books of Breiman and Feller, 
in which the series expansions are replaced by the basic  
convexity inequality 
%$1+x\leq\exp(x)$ for $x\in\R$.
\begin{equation}
  \label{bce}
  \forall t\in\R\qquad
\exp(t)\geq
1+t
\,.
%1+x\leq\exp(x)\,.
\end{equation}
%convexity inequality $1+x\leq e^x$ for $x\in\R$.
%This is absolutely fine, at least in my country, 
This is absolutely fine, at least in France,
where 
%usually 
the exponential function 
is initially introduced to the students as 
the only function which is equal to its derivative and whose 
value at $0$ is~$1$;
inequality~\eqref{bce} can be easily derived from the identity 
%$\exp(t)-1\,=\,\int_0^t\exp(s)\,ds$.
\begin{equation}
  \label{ice}
  \forall t\in\R\qquad
%\exp(t)-1\,=\,\int_0^te^s\,ds\,.
\exp(t)-1\,=\,\int_0^t\exp(s)\,ds\,.
\end{equation}
To sum up, %the prerequisites to read 
the prerequisites necessary to confidently tackle the proof 
that follows
%of the theorem 
are basic knowledge of integral calculus and the 
following elementary inequalities:
\begin{gather}
  \label{ine1}
  \forall t\in[-1,1[\qquad 
\frac{1}{1-t} \,\geq\, 1+t \,,\\
  %\forall x\in]-1/2,1/2[\qquad 
  \label{ine2}
  \forall t\in[0,1/2]\qquad 
\frac{1}{1-t} \,\leq\, 1+2t \,.
\end{gather}
These inequalities are straightforward. 
We simply multiply by $1-t$, develop and remove $1$ from each side. 
The first one reduces to 
$0\geq -t^2$.
%$1\leq 1-t^2$, which holds for $0\leq t<1$.
The second one reduces to 
$0\leq t(1-2t)$, which holds for $0\leq t\leq 1/2$.

With this background knowledge, we will prove completely
theorem~\ref{dml}.
%The proof of de Moivre theorem is done 
We start with de Moivre theorem
in section~\ref{sdm}, which is substantially simpler than the 
general case $p\in]0,1[$.
For an educational presentation in class, our recommendation
is to explain the proof of the upper bound for a symmetric interval.
This corresponds to the first part of 
section~\ref{sdm}, 
that is approximately two pages, knowing that 
%we have written a precise proof with all the intermediate steps.
we have written a detailed proof with all the intermediate steps.
This should fit into a one-hour lecture. 
If one wishes 
to go for the full proof, then one could add the proof of the 
lower bound, which nevertheless has the same flavor as the 
upper bound. However, the lower bound is slightly more difficult,
the inequalities are less straightforward than for the upper bound.
Within a two hours lecture, I could quietly present the theorem and its history,
state the general version of the central limit theorem and perform the detailed proof 
of the upper bound in the symmetric case, and the students 
(at least some of them) looked quite happy!

Modern expositions of the de Moivre--Laplace theorem usually makes 
appeal to the famous Stirling formula (this is the case of both 
Breiman and Feller \cite{B,F}).
If we incorporate Stirling's formula in our initial luggage, we 
arrive faster at the desired conclusion. However, Stirling's 
formula is not necessary at all.
In fact, 
a simpler result, 
namely the formula for the Wallis integral, is more than enough 
%amply enough 
to complete our proof.
As a by product, we
obtain the following non-asymptotic inequality linking 
the binomial and the Gaussian distributions:
\begin{comment}
  \label{gluvex}
  %\forall n\ge\quad \forall n\geq 2x^2\quad%\\
  \forall x>0\quad
\forall n\geq \max(2x^2,1)\qquad%\\
%  P_{2n}(x)
\Big|
P\Big(\big|S_{2n}-{n}\big|<x\sqrt{\frac{n}{2}}\Big)
-
   \int_{-x}^{ x }
   e^{- \tfrac{s^2}{2} }
   %\exp\Big(- \frac{s^2}{2} \Big)\,
   \frac{ds}{
             \sqrt{{2}{\pi}}
   }
\Big|\,\leq\,
   %\times \exp\Big( -\frac{x}{\sqrt{2n}}
           %e^{ \tfrac{x^3+2}{2\sqrt{n}}}-1
           e^{ \tfrac{x^3+2}{\sqrt{n}}}-1
  %{ -\frac{x^3}{\sqrt{2n}}} 
           \,.
\end{comment}
\begin{multline}
  \label{gluvex}
  %\forall n\ge\quad \forall n\geq 2x^2\quad%\\
  \forall x>0\quad
\forall n\geq \max(4x^2,2)\qquad%\\
\cr
%  P_{2n}(x)
\Big|
P\Big(\big|S_{n}-\frac{n}{2}\big|<x{\frac{\sqrt n}{2}}\Big)
-
   \int_{-x}^{ x }
   e^{- \tfrac{s^2}{2} }
   %\exp\Big(- \frac{s^2}{2} \Big)\,
   \frac{ds}{
             \sqrt{{2}{\pi}}
   }
\Big|\,\leq\,
   %\times \exp\Big( -\frac{x}{\sqrt{2n}}
           %e^{ \tfrac{x^3+2}{2\sqrt{n}}}-1
           e^{ \tfrac{4x^3+8}{\sqrt{n}}}-1
  %{ -\frac{x^3}{\sqrt{2n}}} 
           \,.
\end{multline}
\begin{comment}
  \label{gluvex}
  \forall n\geq 1\quad
\forall n\geq 2x^2\qquad
%  P_{2n}(x)
\Bigg|
P\Big(\big|S_{n}-\frac{n}2\big|<x\frac{\sqrt{n}}{2}\Big)
-
   \int_{-x}^{ x }
   e^{- \tfrac{s^2}{2} }
   %\exp\Big(- \frac{s^2}{2} \Big)\,
   \frac{ds}{
             \sqrt{{2}{\pi}}
   }
\Bigg|\,\leq\,
   %\times \exp\Big( -\frac{x}{\sqrt{2n}}
           e^{ \tfrac{x^3+2}{\sqrt{n}}}-1
           \,.
\end{comment}
\begin{comment}
  \label{gluvex}
  \forall n\geq 1\quad
\forall n\geq 2x^2\qquad\\
%  P_{2n}(x)
\Bigg|
P\Big(\big|S_{n}-\frac{n}2\big|<x\frac{\sqrt{n}}{2}\Big)
-
   \int_{-x}^{ x }
   e^{- \tfrac{s^2}{2} }
   %\exp\Big(- \frac{s^2}{2} \Big)\,
   \frac{ds}{
             \sqrt{{2}{\pi}}
   }
\Bigg|\,\leq\,
   %\times \exp\Big( -\frac{x}{\sqrt{2n}}
           e^{ \tfrac{x^3+2}{\sqrt{n}}}-1
           \,.
\end{comment}
%Of course, this inequality is far from optimal, 
%and there exist much finer inequalities of this type in the 
%literature.
The literature abounds with more sophisticated 
inequalities linking the binomial and the Gaussian distributions, 
see for instance \cite{AG,HM,BA}.
One of the most powerful is Tusnady's inequality 
(see \cite{OU} and the references therein).
Hipp and Mattner \cite{HM} provide lower and upper bounds for 
the difference of the 
distribution functions of the symmetric binomial and the Gaussian law 
of the form $c/\sqrt{n}$, with optimal values for the constant $c$.
%Hipp and Mattner \cite{HM} provide lower and upper bounds for 
%the difference of the 
%distribution funcions of the symmetric binomial and the Gaussian law 
%of the form $c/\sqrt{n}$, with optimal values for the constant $c$.
Yet the advantage of the inequality~\eqref{gluvex} lies 
in the simplicity of its proof.

%The interest of the inequalities~\eqref{nakwex},
%~\eqref{uuvex} 
%and~\eqref{luvex} lies in the simplicity of their proofs.

The proof of the general case $p\in]0,1[$ is more involved, it 
is done in section~\ref{dmlt}. 
We would not recommend presenting it in the classroom. 
We include it here, because it 
is simpler than the proofs we found in the existing literature: it does
not require a series expansion, and it relies on the same elementary
inequalities than we used for de Moivre theorem.
%So, even it is not suited for the first class presentation of the de Moivre--Laplace
%theorem, it furnishes an excellent topic for an homewo rk assignment 
However, the use of these inequalities is more tricky, 
because in the symmetric case there is a magical cancellation
which leads to a considerable simplification of the formulas.
Even if this material does not seem suitable for a first presentation 
in class, the proof of the general case would provide an excellent topic 
for a homework assignment, 
allowing students to deepen and generalize the argument of the symmetric case.

\begin{comment}
Finally, we note that
%Furthermore, 
we prove a non-asymptotic inequality linking 
the binomial and the Gaussian distributions.
Of course, the one we obtain here is far from optimal, 
and there exist much finer inequalities of this type in the 
literature: see for instance ?
Yet the advantage of the inequalities presented here lies 
in the simplicity of their proofs.
\end{comment}

%We deal with the symmetric case in
\begin{comment}
So, the purpose of this text is to present an elementary proof of the 
De Moivre--Laplace theorem, which does not make appeal to uniform expansion 
nor to series expansion, and which is accessible to an undergradute student 
once he is familiar with the exponential function and 
the basics of calculus. 
The proof follow the same strategy as the proofs presented
by Breiman and Feller, but the series expansions are replace by adequate 
simpler inequalities. 
\end{comment}
%Although this is a little contribution, 
%we still hope that 
%whose only necessary prerequisites are 
%the rigorous definition of the 
%Based on the proofs presented by Breiman and Feller, 

%\section{The symmetric case.}
%\newpage 
\section{de Moivre theorem.}
\label{sdm}
We consider here the symmetric case $p=1/2$, and also
the case of a symmetric interval. 
Let $x>0$ and let us set 
\begin{equation}
  \forall n\geq 1\qquad 
  P_n(x)\,=\, P\Big(\big|S_n-\frac{n}2\big|<x{\frac{\sqrt n}2}\Big)\,.
\end{equation}
We deal first with the case of even integers $n$, 
which makes the formulas better looking. 
We fix $n\geq 1$ and we start by 
writing down the expression of $P_{2n}(x)$:
\begin{comment}
  \label{dex}
  P_{2n}(x)
  %\,=\, P\Big(\big|S_{2n}-{n}\big|<x\sqrt{\frac{n}2}\Big)\cr
           \,=\, 
           %\sum_{\substack{k\in\zun\\ |k-n|< x\sqrt{\frac{n}2}}}
           %\sum_{\substack{k\in\zdn\\ |k-n|< x\sqrt{{n}/2}}}
           \kern-5pt
           \sum_{k:|k-n|< x\sqrt{{n}/2}}
           \kern-5pt
           P\big(S_{2n}=k\big)
           \,=\, 
           \kern-5pt
           %\sum_{\substack{k\in\zun\\ |k-n|< x\sqrt{\frac{n}2}}}
           %\sum_{\substack{k\in\zdn\\ |k-n|< x\sqrt{{n}/2}}}
           \sum_{k:|k-n|< x\sqrt{{n}/2}}
           \binom{2n}{k}\frac{1}{2^{2n}}
  \,.
\end{comment}
\begin{align}
  \label{dex}
  P_{2n}(x)&\,=\, P\Big(\big|S_{2n}-{n}\big|<x\sqrt{\frac{n}2}\Big)\cr
           &\,=\, 
           %\sum_{\substack{k\in\zun\\ |k-n|< x\sqrt{\frac{n}2}}}
           \sum_{\substack{k\in\zdn\\ |k-n|< x\sqrt{{n}/2}}}
           %\sum_{|k-n|< x\sqrt{{n}/2}}
           P\big(S_{2n}=k\big)
           \,=\, 
           %\sum_{\substack{k\in\zun\\ |k-n|< x\sqrt{\frac{n}2}}}
           \sum_{\substack{k\in\zdn\\ |k-n|< x\sqrt{{n}/2}}}
           %\sum_{|k-n|< x\sqrt{{n}/2}}
           \binom{2n}{k}\frac{1}{2^{2n}}
  \,.
\end{align}
We make the change of variable $j=k-n$, so that~\eqref{dex} becomes
\begin{equation}
  \label{rex}
  P_{2n}(x)
           \,=\, 
           %\sum_{\substack{k\in\zun\\ |k-n|< x\sqrt{\frac{n}2}}}
           \sum_{\substack{j\in\zmn\\ |j|< x\sqrt{{n}/2}}}
           \binom{2n}{n+j}\frac{1}{2^{2n}}
  \,.
\end{equation}
We take apart the term for $j=0$, and we use the symmetry of the 
binomial coefficient to rewrite~\eqref{rex} as
\begin{equation}
  \label{vex}
  P_{2n}(x)
           \,=\, 
           %\sum_{\substack{k\in\zun\\ |k-n|< x\sqrt{\frac{n}2}}}
           \frac{1}{2^{2n}}
           \binom{2n}{n}
           +
          2 \sum_{\substack{j\in\unn\\ j< x\sqrt{{n}/2}}}
           \binom{2n}{n+j}\frac{1}{2^{2n}}
  \,.
\end{equation}
This formula is the starting point of the computations.
For $j\geq 1$, we have 
\begin{align}
  \label{fact1}
  (n+j)!&\,=\,n!\,(n+1)\cdots(n+j)\,,\\ \label{fact2}
  (n-j)!&\,=\,\frac{n!}{(n-j+1)\cdots n}\,,
\end{align}
As in the classical arguments of Breiman and Feller, 
we use~\eqref{fact1} and~\eqref{fact2} 
to rewrite the binomial coefficient 
           $\binom{2n}{n+j}$
with the help of the central binomial coefficient
           $\binom{2n}{n}$
           as follows:
\begin{equation}
  \label{cent}
           \binom{2n}{n+j}\,=\,
           \frac{(2n)!}{{(n+j)!}{(n-j)!}}
           \,=\,
           %\frac{(2n)!}{n! n!}
           \binom{2n}{n}
  \frac{(n-j+1)\cdots n}
  {(n+1)\cdots(n+j)}
  \,.
\end{equation}
We introduce next a notation for the product appearing in~\eqref{cent},
\begin{equation}
  \label{pent}
  \forall j\geq 1\qquad
  \pi(j,n)
           \,=\,
  \frac{(n-j+1)\cdots n}
  {(n+1)\cdots(n+j)}
           \,=\,
           \prod_{1\leq\ell\leq j}
  \frac{n+\ell-j}{n+\ell}
  \,.
\end{equation}
We set also $\pi(0,n)=1$.
%With this new notation, the formula~\eqref{vex} can be rewritten as
With this new notation, the formula~\eqref{vex} becomes
\begin{equation}
  \label{bvex}
  P_{2n}(x)
           \,=\, 
           %\sum_{\substack{k\in\zun\\ |k-n|< x\sqrt{\frac{n}2}}}
           %\frac{2}{2^{2n}} \binom{2n}{n}
           P(S_{2n}=n)
           \Bigg(
           1%  \frac12
           +
           2\sum_{\substack{j\in\unn\\ j< x\sqrt{{n}/2}}}
           \pi(j,n)
           \Bigg)
  \,.
\end{equation}
%\subsection{lower bound of $\pi(j,n)$.}
%\subsection{Upper bound of $\pi(j,n)$.}
\subsection{Upper bound.}
%Let us fix $j\geq 1$. To prepare for the lower bound, we rewrite 
Let us fix $j\in\unn$. To prepare for the upper bound, we rewrite 
$\pi(j,n)$ as follows:
\begin{equation}
  \label{manop}
  \pi(j,n) \,=\, \prod_{1\leq\ell\leq j}
  \frac{n+\ell-j}{n+\ell}
   \,=\, \prod_{1\leq\ell\leq j}
   \Big(1-
   \frac{j}{n+\ell}
     \Big)
  \,,
\end{equation}
and we use inequality~\eqref{bce} 
to deduce from~\eqref{manop} that
\begin{equation}
  \label{manor}
  \pi(j,n) \,\leq\, 
  %\prod_{1\leq\ell\leq j}
  %\exp\Big(-
  %\frac{j}{n+\ell}
    %\Big)
   %\,=\, 
  \exp\Big(
  - \sum_{1\leq\ell\leq j}
   \frac{j}{n+\ell}
     \Big)
   \,\leq\, 
  \exp\Big(
  - 
   \frac{j^2}{n+j}
     \Big)
  \,.
\end{equation}
Next, we use the inequality~\eqref{ine1}
\begin{comment}
  \forall t\in\Big[0,\frac{1}{2}\Big]\qquad 
  \frac{1}{1-t}\,\leq\, 1+2t
\end{comment}
to obtain 
\begin{equation}
  \label{mbnos}
   \frac{j^2}{n+j}
   \,=\,
   \frac{j^2}{n}
   \frac{1}{ \displaystyle 1+\frac{j}{n}
   }
   \,\geq\,
   \frac{j^2}{n}
   \Big( 1-\frac{j}{n}\Big)
  \,.
\end{equation}
Substituting successively~\eqref{mbnos} in~\eqref{manor}
and~\eqref{manop}, we conclude that
\begin{equation}
  \label{hinos}
  \pi(j,n) 
  \,\leq \, 
   \exp\Big(-
   \frac{j^2}{n}
   +\frac{j^3}{n^2}
   \Big)
  \,.
\end{equation}
\begin{comment}
Contrary to the previous case, this inequality is as good as the one 
obtained by 
Feller 
(see \cite{F}, page 180, between formulas~$(2.3)$ and $(2.4)$).
\end{comment}
%\subsection{Upper bound of $P_{2n}(x)$.}
With the upper bound of $\pi(j,n)$ in hand, we compute 
%With this inequalitye upper bound of $\pi(j,n)$ in hand, 
%We compute next
an upper bound on 
$P_{2n}(x)$. Let us fix $x>0$ and let $n$ be an integer 
such that $n\geq x^2$. This ensures that 
\begin{equation}
  \label{fnst}
           1\leq j< x\sqrt{\frac{n}2}\quad
           \Longrightarrow
           \quad
           {j}\leq n\,,\quad
         %  \frac{j}{n}\leq\frac12\,,\quad
           \frac{j^3}{n^2}\leq\frac{x^3}{2\sqrt{2n}}\,.
\end{equation}
%In particular, for these values of $j$, the inequality~\eqref{hinos}
For these values of $j$, the inequality~\eqref{hinos}
holds, and we have, thanks to~\eqref{fnst},
\begin{equation}
  \label{rvinos}
           1\leq j< x\sqrt{\frac{n}2}\quad
           \Longrightarrow
           \quad 
  \pi(j,n) 
  \,\leq \, 
   \exp\Big(-
   \frac{j^2}{n}
 +\frac{x^3}{2\sqrt{2n}}
   \Big)
  \,.
\end{equation}
%we can substitute inequality~\eqref{linos}
%In particular, we can substitute inequality~\eqref{linos}
Substituting inequality~\eqref{rvinos}
into the sum appearing in~\eqref{bvex}, we get
\begin{equation}
  \label{tvex}
  P_{2n}(x)
           \,\leq\, 
           %\sum_{\substack{k\in\zun\\ |k-n|< x\sqrt{\frac{n}2}}}
           P(S_{2n}=n)
           \Bigg(
           1%  \frac12
           +
            2 e^{ \frac{x^3}{2\sqrt{2n}}}
           \sum_{1\leq j< x\sqrt{{n}/2}}
           %\sum_{\substack{j\in\unn\\ j< x\sqrt{{n}/2}}}
   \exp\Big(-
   \frac{j^2}{n}
 %+\frac{x^3}{2\sqrt{2n}}
   \Big)
           \Bigg)
  \,.
\end{equation}
\begin{comment}
  \label{tvex}
  P_{2n}(x)
           &\,\leq\, 
           %\sum_{\substack{k\in\zun\\ |k-n|< x\sqrt{\frac{n}2}}}
           2P(S_{2n}=n)
           \Bigg(
             \frac12
           +
           \sum_{\substack{j\in\unn\\ j< x\sqrt{{n}/2}}}
   \exp\Big(-
   \frac{j^2}{n}
 +\frac{x^3}{2\sqrt{2n}}
   \Big)
           \Bigg)\cr
           &         \,\geq\, 
           %\sum_{\substack{k\in\zun\\ |k-n|< x\sqrt{\frac{n}2}}}
           2P(S_{2n}=n)
           \Bigg(
   %\exp\Big( -\frac{x}{\sqrt{2n}} \Big)
             e^{ -\frac{x}{\sqrt{2n}}}
           \sum_{\substack{j\in\zun\\ j< x\sqrt{{n}/2}}}
   \exp\Big(-
   \frac{j^2}{n}
   \Big)
             -\frac12
           \Bigg)
  \,.
\end{comment}
%So, why did we bother to rederive 
%inequality~\eqref{tmyin} from scratch? 
%It should be noted that the proof of Keller and Kindler
%is quite complex and difficult (although not as mysterious
%as Talagrand's proof), and it involves also other
%deep results, typically hypercontractivity estimates.
Since the function $\exp(-t)$ is decreasing on $[0,+\infty[$, then
\begin{equation}
  \label{dosi}
  \forall j\geq 1\qquad 
   \exp\Big(- \frac{j^2}{n} \Big)\,\leq\,
   \int_{j-1}^{j}
   e^{- \tfrac{t^2}{n} }
  % \exp\Big(- \frac{t^2}{n} \Big)\,
   dt\,.
\end{equation}
Substituting inequality~\eqref{dosi}
into the sum appearing in~\eqref{tvex}, we get
\begin{align}
  \label{fwex}
  P_{2n}(x)
           &\,\leq\, 
           %\sum_{\substack{k\in\zun\\ |k-n|< x\sqrt{\frac{n}2}}}
           P(S_{2n}=n)
           \Bigg(
           1%  \frac12
           +
           2  e^{ \frac{x^3}{2\sqrt{2n}}}
           %\sum_{\substack{j\in\unn\\ j< x\sqrt{{n}/2}}}
           \sum_{1\leq  j< x\sqrt{{n}/2}}
   \int_{j-1}^{j}
   e^{- \tfrac{t^2}{n} }
  % \exp\Big(- \frac{t^2}{n} \Big)\,
   dt
   %\times \exp\Big( -\frac{x}{\sqrt{2n}}
           \Bigg)
           \cr
           &\,\leq\, 
           %\sum_{\substack{k\in\zun\\ |k-n|< x\sqrt{\frac{n}2}}}
           P(S_{2n}=n)
           \Bigg(
           1%  \frac12
           +
           2  e^{ \frac{x^3}{2\sqrt{2n}}}
   \int_{0}^{ x\sqrt{{n}/2} }
   e^{- \tfrac{t^2}{n} }
  % \exp\Big(- \frac{t^2}{n} \Big)\,
   dt
   %\times \exp\Big( -\frac{x}{\sqrt{2n}}
           \Bigg)
  \,.
\end{align}
We make the change of variable 
$s { \sqrt{{n}/2} }=t$ in~\eqref{fwex} and the inequality becomes
\begin{equation}
  \label{wvex}
  P_{2n}(x)
           \,\leq\, 
           %\sum_{\substack{k\in\zun\\ |k-n|< x\sqrt{\frac{n}2}}}
           %\sqrt{2n}P(S_{2n}=n)
           P(S_{2n}=n)
           \Bigg(
          1%   \frac12
           +
             %\sqrt{\frac{n}{2}}
             \sqrt{2{n}}\,
             e^{ \frac{x^3}{2\sqrt{2n}}}
   \int_{0}^{ x }
   e^{- \tfrac{s^2}{2} }
  % \exp\Big(- \frac{s^2}{2} \Big)\,
   ds
   %\times \exp\Big( -\frac{x}{\sqrt{2n}}
           \Bigg)
  \,.
\end{equation}
We finally factorize $\sqrt{n}$, we use the symmetry of 
the integrand, and we 
introduce the convenient normalizing factor $\sqrt{2\pi}$:
\begin{equation}
  \label{uwex}
  P_{2n}(x)
           \,\leq\, 
           %\sum_{\substack{k\in\zun\\ |k-n|< x\sqrt{\frac{n}2}}}
           %\sqrt{2n}P(S_{2n}=n)
             \sqrt{{n}{\pi}}
           P(S_{2n}=n)
           \Bigg(
             \frac{1}{ \sqrt{{n}{\pi}} }
           +
             e^{ \frac{x^3}{2\sqrt{2n}}}
   \int_{-x}^{ x }
   %\exp\Big(- \frac{s^2}{2} \Big)\,
   e^{- \tfrac{s^2}{2} }
   \frac{ds}{
             \sqrt{{2}{\pi}}
   }
   %\times \exp\Big( -\frac{x}{\sqrt{2n}}
           \Bigg)
  \,.
\end{equation}
%\subsection{Lower bound of $\pi(j,n)$.}
\subsection{Lower bound.}
%Let us fix $j\geq 1$. To prepare for the lower bound, we rewrite 
Let us fix $j\in\unns$. To prepare for the lower bound, we rewrite 
$\pi(j,n)$ as follows:
\begin{equation}
  \label{minop}
  \pi(j,n) \,=\, \prod_{1\leq\ell\leq j}
  \frac{n+\ell-j}{n+\ell}
   \,=\, \prod_{1\leq\ell\leq j}
   \frac{1}{\displaystyle 1+\frac{j}{n+\ell-j}}
  \,.
\end{equation}
Notice that
\begin{equation}
  \label{nott}
  \forall \ell\in
{\{\,1,\dots,j\,\}}\qquad
   \frac{j}{n+\ell-j}\,\leq\,
   \frac{j}{n+1-j}\,\leq\,
   \frac{n/2}{n/2+1}\,<\,1
   \,.
\end{equation}
The inequality~\eqref{bce} readily implies that
\begin{equation}
  \label{lowi}
  \forall t\in[0,1]\qquad 
  \frac{1}{1+t}\,\geq\, \exp(-t) \,,
\end{equation}
and we use this inequality 
to 
deduce from~\eqref{nott} and~\eqref{minop} that
\begin{equation}
  \label{minos}
  \pi(j,n) 
  \,\geq \, \prod_{1\leq\ell\leq j}
   \exp\Big(-\frac{j}{n+\ell-j}\Big)
   \,=\,\exp\Big(-
   \sum_{1\leq\ell\leq j}
   \frac{j}{n+\ell-j}\Big)
  \,.
\end{equation}
Next, 
we bound crudely from above the argument of the exponential :
\begin{equation}
  \label{abnos}
   \sum_{1\leq\ell\leq j}
   \frac{j}{n+\ell-j}
   \,\leq\,
   \frac{j^2}{n-j}
   \,,
\end{equation}
and we use the inequality~\eqref{ine2}
\begin{comment}
  \forall t\in\Big[0,\frac{1}{2}\Big]\qquad 
  \frac{1}{1-t}\,\leq\, 1+2t
\end{comment}
to obtain 
\begin{equation}
  \label{bbnos}
   \frac{j^2}{n-j}
   \,=\,
   \frac{j^2}{n}
   \frac{1}{ \displaystyle 1-\frac{j}{n}
   }
   \,\leq\,
   \frac{j^2}{n}
   \Big( 1+2\frac{j}{n}\Big)
  \,.
\end{equation}
Substituting successively~\eqref{bbnos} in~\eqref{abnos}
and~\eqref{minos}, we conclude that
\begin{equation}
  \label{linos}
  \pi(j,n) 
  \,\geq \, 
   \exp\Big(-
   \frac{j^2}{n}
   -2\frac{j^3}{n^2}
   \Big)
  \,.
\end{equation}
Notice that Feller obtained the better inequality where the $2$ in~\eqref{linos}
is replaced by~$1$ (see \cite{F}, page 180, 
between formulas~$(2.3)$ and $(2.4)$).
However, his computation relies on a series 
expansion and a comparison with  an adequate geometric series.
For the upper bound, we obtained the inequality~\eqref{hinos},
which is the same as the one of Feller.
%without the factor~$2$.
%\subsection{lower bound of $P_{2n}(x)$.}
%\subsection{Lower bound of $P_{2n}(x)$.}

With the lower bound of $\pi(j,n)$ in hand, we compute 
a lower bound on 
$P_{2n}(x)$. Let us fix $x>0$ and let $n$ be an integer 
such that $n\geq 2x^2$. This ensures that 
\begin{equation}
  \label{enst}
           1\leq j< x\sqrt{\frac{n}2}\quad
           \Longrightarrow
           \quad
           %\frac{j}{n}\leq\frac12\,,\quad
           {j}\leq\frac{n}2\,,\quad
           \frac{j^3}{n^2}\leq\frac{x^3}{2\sqrt{2n}}\,.
\end{equation}
%In particular, for these values of $j$, the inequality~\eqref{linos}
\begin{comment}
For these values of $j$, the inequality~\eqref{linos}
holds, and we have, thanks to~\eqref{enst},
\begin{equation}
  \label{yvinos}
           1\leq j< x\sqrt{\frac{n}2}\quad
           \Longrightarrow
           \quad 
  \pi(j,n) 
  \,\geq \, 
   \exp\Big(-
   \frac{j^2}{n}
 -\frac{x^3}{\sqrt{2n}}
   \Big)
  \,.
\end{equation}
%we can substitute inequality~\eqref{linos}
%In particular, we can substitute inequality~\eqref{linos}
\end{comment}
%Substituting inequality~\eqref{yvinos}
%For these values of $j$, 
Using the inequalities~\eqref{linos} and~\eqref{enst}
into the sum appearing in~\eqref{bvex}, we get
\begin{align}
  \label{hvex}
  P_{2n}(x)
           &\,\geq\, 
           %\sum_{\substack{k\in\zun\\ |k-n|< x\sqrt{\frac{n}2}}}
           P(S_{2n}=n)
           \Bigg(
            1% \frac12
           +
           %2\sum_{\substack{j\in\unn\\ j< x\sqrt{{n}/2}}}
           2\sum_{1\leq j< x\sqrt{{n}/2}}
   \exp\Big(-
   \frac{j^2}{n}
 -\frac{x^3}{\sqrt{2n}}
   \Big)
           \Bigg)\cr
           &         \,\geq\, 
           %\sum_{\substack{k\in\zun\\ |k-n|< x\sqrt{\frac{n}2}}}
           P(S_{2n}=n)
           \Bigg(
   %\exp\Big( -\frac{x}{\sqrt{2n}} \Big)
             2e^{ -\frac{x^3}{\sqrt{2n}}}
           %\sum_{\substack{j\in\zun\\ j< x\sqrt{{n}/2}}}
           \sum_{0\leq j< x\sqrt{{n}/2}}
   \exp\Big(-
   \frac{j^2}{n}
   \Big)
            -1 % -\frac12
           \Bigg)
  \,.
\end{align}
%So, why did we bother to rederive 
%inequality~\eqref{tmyin} from scratch? 
%It should be noted that the proof of Keller and Kindler
%is quite complex and difficult (although not as mysterious
%as Talagrand's proof), and it involves also other
%deep results, typically hypercontractivity estimates.
Since the function $\exp(-t)$ is decreasing on $[0,+\infty[$, then
\begin{equation}
  \label{cosi}
  \forall j\geq 0\qquad 
   \exp\Big(- \frac{j^2}{n} \Big)\,\geq\,
   \int_{j}^{j+1}
   e^{- \tfrac{t^2}{n} }
   %\exp\Big(- \frac{t^2}{n} \Big)\,
   dt\,.
\end{equation}
Substituting inequality~\eqref{cosi}
into the sum appearing in~\eqref{hvex}, we get
\begin{align}
  \label{fvex}
  P_{2n}(x)
           &\,\geq\, 
           %\sum_{\substack{k\in\zun\\ |k-n|< x\sqrt{\frac{n}2}}}
           P(S_{2n}=n)
           \Bigg(
            2 e^{ -\frac{x^3}{\sqrt{2n}}}
           %\sum_{\substack{j\in\zun\\ j< x\sqrt{{n}/2}}}
           \sum_{0\leq j< x\sqrt{{n}/2}}
   \int_{j}^{j+1}
   e^{- \tfrac{t^2}{n} }
  % \exp\Big(- \frac{t^2}{n} \Big)\,
   dt
   %\times \exp\Big( -\frac{x}{\sqrt{2n}}
           -1 %  -\frac12
           \Bigg)
           \cr
           &\,\geq\, 
           %\sum_{\substack{k\in\zun\\ |k-n|< x\sqrt{\frac{n}2}}}
           P(S_{2n}=n)
           \Bigg(
             2e^{ -\frac{x^3}{\sqrt{2n}}}
   \int_{0}^{ x\sqrt{{n}/2} }
   e^{- \tfrac{t^2}{n} }
  % \exp\Big(- \frac{t^2}{n} \Big)\,
   dt
   %\times \exp\Big( -\frac{x}{\sqrt{2n}}
             -1 %-\frac12
           \Bigg)\cr
           &\,\geq\, 
           %\sum_{\substack{k\in\zun\\ |k-n|< x\sqrt{\frac{n}2}}}
           %\sqrt{2n}P(S_{2n}=n)
             \sqrt{{n}{\pi}}
           P(S_{2n}=n)
           \Bigg(
             e^{ -\frac{x^3}{\sqrt{2n}}}
   \int_{-x}^{ x }
   e^{- \tfrac{s^2}{2} }
   %\exp\Big(- \frac{s^2}{2} \Big)\,
   \frac{ds}{
             \sqrt{{2}{\pi}}
   }
   %\times \exp\Big( -\frac{x}{\sqrt{2n}}
             -
             \frac{1}{
             \sqrt{{n}{\pi}}
             }
           \Bigg)
  \,,
\end{align}
where we have made the same final steps as in the upper bound 
to get the last formula, namely 
the change of variable 
$s { \sqrt{{n}/2} }=t$, the use of symmetry of the integrand and 
factorization of $\sqrt{n\pi}$. 
\subsection{Conclusion.}
It follows from Stirling's formula that 
\begin{equation}
  \label{fstir}  
  \lim_{n\to\infty}
             \sqrt{{n}{\pi}}
           P(S_{2n}=n)\,=\,1\,.
\end{equation}
Using~\eqref{fstir} and passing to the limit in 
the inequalities~\eqref{uwex} and~\eqref{fvex}, we conclude that
\begin{equation}
  \label{concl}
  \lim_{n\to\infty}
  P_{2n}(x)
           \,=\, 
   \int_{-x}^{ x }
   e^{- \tfrac{s^2}{2} }
   %\exp\Big(- \frac{s^2}{2} \Big)\,
   \frac{ds}{
             \sqrt{{2}{\pi}}
   }
  \,.
\end{equation}
\subsection{The case of odd integers.} So far, we have dealt 
only with the subsequence of the even integers.
For the odd integers, we could proceed as Breiman to show that 
the infimum limit and the supremum limit of 
the sequence $(P_n(x))_{n\geq 0}$
coincide along the even and odd integers.
             %$\sqrt{{n}{\pi}} P(S_{2n}=n)$ 
             %$\sqrt{{n}{\pi}} P(S_{2n+1}=n)$ 
Or we could simply say that 
the previous inequalities can be adapted with minor modifications 
to handle the case of odd integers. However, this would be a bit unfair,
because there do exist differences between the two cases, 
and the work has to be done if one really wishes to obtain the 
correct inequalities. So we present next an accelerated version 
of the derivation of these inequalities, which rests nevertheless 
on the inequalities derived in the even case.
%the sequence 
%$(P_{2n+1}(x))_{n\geq 0}$ (CHECK or DO). 

Let us fix $n\geq 0$. We have
\begin{multline}
  \label{odex}
  P_{2n+1}(x)\,=\, P\Big(\big|S_{2n+1}-\frac{2n+1}{2}\big|<
  \frac{x}{2}\sqrt{2n+1}
\Big)\cr
           \,=\, 
           %\sum_{\substack{k\in\zun\\ |k-n|< x\sqrt{\frac{n}2}}}
           \sum_{\substack{k\in\zdnp\\ |k-\frac{2n+1}2|< 
  \frac{x}{2}\sqrt{2n+1}
           }}
           P\big(S_{2n+1}=k\big)
           \,=\, 
           %\sum_{\substack{k\in\zun\\ |k-n|< x\sqrt{\frac{n}2}}}
           \sum_{\substack{k\in\zdnp\\ |k-n-\frac12|< 
  \frac{x}{2}\sqrt{2n+1}
           }}
           \binom{2n+1}{k}\frac{1}{2^{2n+1}}\cr
           \,=\, 
           %\sum_{\substack{k\in\zun\\ |k-n|< x\sqrt{\frac{n}2}}}
           \sum_{\substack{k\in\zun\\ n+\frac12-k< 
  \frac{x}{2}\sqrt{2n+1}
           }}
           \binom{2n+1}{k}\frac{1}{2^{2n+1}}
           \,+\, 
           %\sum_{\substack{k\in\zun\\ |k-n|< x\sqrt{\frac{n}2}}}
           \sum_{\substack{k\in \{\,n+1,\dots,2n+1\,\}
               \\ k-n-\frac12< 
  \frac{x}{2}\sqrt{2n+1}
           }}
           \binom{2n+1}{k}\frac{1}{2^{2n+1}}
  \,.
\end{multline}
We perform the change of variables $j=n-k$ in the first sum and 
$j=k-n-1$ in the second sum.
Using the symmetry of the 
binomial coefficient, we see that the two sums are in fact equal and we get 
\begin{equation}
  \label{ovex}
  P_{2n+1}(x)
           \,=\, 
           %\sum_{\substack{k\in\zun\\ |k-n|< x\sqrt{\frac{n}2}}}
           2\sum_{\substack{j \in\zun\\ j+\frac12< 
  \frac{x}{2}\sqrt{2n+1}
           }}
           \binom{2n+1}{n-j}\frac{1}{2^{2n+1}}
           %\binom{2n+1}{n+1+j}\frac{1}{2^{2n+1}}
  \,.
\end{equation}
This formula is the analog for the odd case of formula~\eqref{vex}. Notice that 
we did not have to isolate the central term, because when $n$ is even, 
the symmetric binomial distribution
possesses two central terms, which are naturally incorporated in the sum.
We rewrite next the binomial coefficient 
$\smash{\binom{2n+1}{n-j}}$
with the help of the central binomial coefficient
%in terms of the central binomial coefficient
$\smash{\binom{2n+1}{n}}$
           as follows:
\begin{equation}
  \label{ocent}
           \binom{2n+1}{n-j}\,=\,
           \frac{(2n+1)!} {(n-j)! {(n+1+j)!} }
           \,=\,
           %\frac{(2n)!}{n! n!}
           \binom{2n+1}{n}
  \frac{(n-j+1)\cdots n}
  {(n+2)\cdots(n+j+1)}
  \,.
\end{equation}
We define 
\begin{equation}
  \label{opent}
  \forall j\geq 0\qquad
  \tpi(j,n)
           \,=\,
  \frac{(n-j+1)\cdots n}
  {(n+2)\cdots(n+j+1)}
           \,=\,
           \prod_{1\leq\ell\leq j}
  \frac{n+\ell-j}{n+\ell+1 }
  \,,
\end{equation}
and the formula~\eqref{ovex} becomes
\begin{equation}
  \label{obvex}
  P_{2n+1}(x)
           \,=\, 
           %\sum_{\substack{k\in\zun\\ |k-n|< x\sqrt{\frac{n}2}}}
           %\frac{2}{2^{2n}} \binom{2n}{n}
           2P(S_{2n+1}=n)
           \sum_{\substack{j \in\zun\\ j+\frac12< 
  \frac{x}{2}\sqrt{2n+1}
           }}
  \tpi(j,n)
           %\frac{1}{2^{2n+1}}
           %\binom{2n+1}{n+1+j}\frac{1}{2^{2n+1}}
  \,.
\end{equation}
%\subsection{lower bound of $\pi(j,n)$.}
%\subsection{Upper bound of $\pi(j,n)$.}
We could then proceed as in the even case to bound 
$\tpi(j,n)$ from above and from below.
Instead, we notice that
\begin{equation}
  \label{adv}
\forall j\in\zun\qquad 
\tpi(j,n)\,=\,
\frac{n+1}{n+j+1}
\pi(j,n)\,,
\end{equation}
%we use the simple inequalities 
\begin{comment}
  \label{adv}
\forall j\in\zun\qquad 
\frac{n}{n+j+1}
\pi(j,n)\,\leq\,
\tpi(j,n)\,\leq\,
\pi(j,n)\,,
\end{comment}
and we take advantage of the work already done on $\pi(j,n)$: 
the inequalities~\eqref{hinos},\eqref{linos} 
and the identity~\eqref{adv} yield that,
for $j\in\unns$, 
\begin{equation}
  \label{ytnos}
\frac{n+1}{n+j+1}
   \exp\Big(-
   \frac{j^2}{n}
   -2\frac{j^3}{n^2}
   \Big)
  \,\leq \, 
\tpi(j,n) 
  \,\leq \, 
   \exp\Big(-
   \frac{j^2}{n}
   +\frac{j^3}{n^2}
   \Big)
  \,.
\end{equation}
Let us fix $x>0$ and let $n$ be an integer 
such that $n\geq x^2$. 
Setting
\begin{equation}
  \label{xnd}
x_n\,=\, {\frac{x}{2}\sqrt{2n+1}}-\frac12\,,
\end{equation}
we have, thanks to the concavity of the square root function,
\begin{equation}
  \label{concavs}
  %\frac{x}{2}\sqrt{2n+1} -\frac12\,\leq\,
  x_n\,\leq\,
           %\frac{x}{2}\sqrt{2n}\sqrt{1+\frac{1}{2n}} -\frac12\,\leq\,
           x\sqrt{\frac{n}{2}}\Big(1+\frac{1}{4n}\Big)
           -\frac12
           \,\leq\,
           x\sqrt{\frac{n}{2}}
           +
           \frac{x}{\sqrt{32n}}
           -\frac12\,<\,
           x\sqrt{\frac{n}{2}}\,.
\end{equation}
\begin{comment}
  \label{concavs}
  %\frac{x}{2}\sqrt{2n+1} -\frac12\,\leq\,
  x_n\,\leq\,
           \frac{x}{2}\sqrt{2n}\sqrt{1+\frac{1}{2n}}
           -\frac12\,\leq\,
           x\sqrt{\frac{n}{2}}\Big(1+\frac{1}{4n}\Big)
           -\frac12\cr 
           \,\leq\,
           x\sqrt{\frac{n}{2}}
           +
           \frac{x}{\sqrt{32n}}
           -\frac12\,<\,
           x\sqrt{\frac{n}{2}}\,.
\end{comment}
It follows from~\eqref{xnd} and~\eqref{concavs} that
\begin{comment}
whenever 
           $0<j+\frac12< 
  \frac{x}{2}\sqrt{2n+1}$, then
           $0\leq j< x\sqrt{{n}/2}$ and
\end{comment}
\begin{equation}
           0<j+\frac12< 
  \frac{x}{2}\sqrt{2n+1}\quad\Longrightarrow \quad
           0\leq j< x\sqrt{\frac{n}2}\,.
\end{equation}
           %$1\leq j< x\sqrt{{n}/2}$, then
           %$1\leq j< x\sqrt{\frac{n}2}$, then
Thus we can use~\eqref{fnst} to conclude from~\eqref{ytnos} that, 
for the indices $j$ appearing in the sum~\eqref{obvex}, we have
\begin{equation}
  \label{ztnos}
  \frac{n}{n+x\sqrt{n}}
   \exp\Big(-
   \frac{j^2}{n}
   -\frac{x^3}{\sqrt{2n}}
   \Big)
  \,\leq \, 
\tpi(j,n) 
  \,\leq \, 
   \exp\Big(-
   \frac{j^2}{n}
   +\frac{x^3}{2\sqrt{2n}}
   \Big)
  \,.
\end{equation}
Substituting inequalities~\eqref{ztnos}
into the sum~\eqref{obvex}, 
and performing the same steps as for the even case, we get
\begin{equation}
  \label{ofwex}
             e^{- \frac{x^3}{\sqrt{2n}}}
   \int_{0}^{x_n+1}
   e^{- \tfrac{t^2}{n} }
  % \exp\Big(- \frac{t^2}{n} \Big)\,
   dt
   %\times \exp\Big( -\frac{x}{\sqrt{2n}}
           \,\leq\, 
           \frac{ P_{2n+1}(x)}{
           2P(S_{2n+1}=n)
           }
           %\sum_{\substack{k\in\zun\\ |k-n|< x\sqrt{\frac{n}2}}}
           \,\leq\, 
           %\sum_{\substack{k\in\zun\\ |k-n|< x\sqrt{\frac{n}2}}}
             1
           +
             e^{ \frac{x^3}{2\sqrt{2n}}}
   \int_{0}^{x_n}
   e^{- \tfrac{t^2}{n} }
  % \exp\Big(- \frac{t^2}{n} \Big)\,
   dt
   %\times \exp\Big( -\frac{x}{\sqrt{2n}}
  \,.
\end{equation}
%where we have set $x_n= {\frac{x}{2}\sqrt{2n+1}}+\frac12$.
We see from~\eqref{xnd} that 
 $x_n+1\geq x\sqrt{{n}/{2}}$.
%Notice that 
\begin{comment}
 x\sqrt{\frac{n}{2}}\,\leq\,x_n+1\,,\qquad
 x_n
 \,\leq\,
 \sqrt{\frac{n}{2}}
 \Big( x+\frac{x}{4n}-
 \frac{1}{\sqrt{2n}}\Big)
 \,\leq\,
 \sqrt{\frac{n}{2}}x
 \,,
\end{comment}
%where the last inequality holds when $x\leq\sqrt{8n}$.
\begin{comment}
  \label{ofwex}
           2P(S_{2n+1}=n)
           \Bigg(
             e^{- \frac{x^3}{\sqrt{2n}}}
   \int_{0}^{\frac{x}{2}\sqrt{2n+1}+1/2}
   e^{- \tfrac{t^2}{n} }
  % \exp\Big(- \frac{t^2}{n} \Big)\,
   dt
   %\times \exp\Big( -\frac{x}{\sqrt{2n}}
           \Bigg)
           \,\leq\, \cr
  P_{2n+1}(x)
           %\sum_{\substack{k\in\zun\\ |k-n|< x\sqrt{\frac{n}2}}}
           \,\leq\, 
           %\sum_{\substack{k\in\zun\\ |k-n|< x\sqrt{\frac{n}2}}}
           2P(S_{2n+1}=n)
           \Bigg(
             1
           +
             e^{ \frac{x^3}{2\sqrt{2n}}}
   \int_{0}^{\frac{x}{2}\sqrt{2n+1}-1/2}
   e^{- \tfrac{t^2}{n} }
  % \exp\Big(- \frac{t^2}{n} \Big)\,
   dt
   %\times \exp\Big( -\frac{x}{\sqrt{2n}}
           \Bigg)
  \,.
\end{comment}
%Assuming it is the case,
We use also the inequality~\eqref{concavs},
we make the change of variable 
$s { \sqrt{{n}/2} }=t$ in~\eqref{ofwex} and we get
\begin{equation}
  \label{owvex}
             e^{- \frac{x^3}{\sqrt{2n}}}
             \sqrt{\frac{n}{2}}
   \int_{0}^{ x }
   e^{- \tfrac{s^2}{2} }
  % \exp\Big(- \frac{s^2}{2} \Big)\,
   ds
           \,\leq\, 
           \frac{ P_{2n+1}(x)}{
           2P(S_{2n+1}=n)
           }
           %\sum_{\substack{k\in\zun\\ |k-n|< x\sqrt{\frac{n}2}}}
           \,\leq\, 
           %\sum_{\substack{k\in\zun\\ |k-n|< x\sqrt{\frac{n}2}}}
             1
             %\frac12
           +
             \sqrt{\frac{n}{2}}
             e^{ \frac{x^3}{2\sqrt{2n}}}
   \int_{0}^{ x }
   e^{- \tfrac{s^2}{2} }
  % \exp\Big(- \frac{s^2}{2} \Big)\,
   ds
  \,.
\end{equation}%the inequality becomes, for $n\geq 8x^2$,
We finally factorize $\sqrt{n}$, we use the symmetry of 
the integrand, and we 
introduce the convenient normalizing factor $\sqrt{2\pi}$:
\begin{equation}
  \label{oouwvex}
             e^{- \frac{x^3}{\sqrt{2n}}}
   \int_{-x}^{ x }
   %\exp\Big(- \frac{s^2}{2} \Big)\,
   e^{- \tfrac{s^2}{2} }
   \frac{ds}{
             \sqrt{{2}{\pi}}
   }
  % \exp\Big(- \frac{s^2}{2} \Big)\,
   ds
           \,\leq\, 
           \frac{ P_{2n+1}(x)}{
             { \sqrt{{n}{\pi}} }
           P(S_{2n+1}=n)
           }
           %\sum_{\substack{k\in\zun\\ |k-n|< x\sqrt{\frac{n}2}}}
           \,\leq\, 
           %\sum_{\substack{k\in\zun\\ |k-n|< x\sqrt{\frac{n}2}}}
           \frac{2}{\sqrt n}
             %\frac12
           +
             e^{ \frac{x^3}{2\sqrt{2n}}}
   \int_{-x}^{ x }
   %\exp\Big(- \frac{s^2}{2} \Big)\,
   e^{- \tfrac{s^2}{2} }
   \frac{ds}{
             \sqrt{{2}{\pi}}
   }
  \,.
\end{equation}
The inequalities~\eqref{oouwvex} readily imply that
%\subsection{Lower bound of $\pi(j,n)$.}
\begin{equation}
  \label{oconcl}
  \lim_{n\to\infty}
  P_{2n+1}(x)
           \,=\, 
   \int_{-x}^{ x }
   e^{- \tfrac{s^2}{2} }
   %\exp\Big(- \frac{s^2}{2} \Big)\,
   \frac{ds}{
             \sqrt{{2}{\pi}}
   }
  \,,
\end{equation}
and it follows finally from~\eqref{concl} and~\eqref{oconcl} that
\begin{equation}
  \lim_{n\to\infty}
  P\Big(\big|S_n-\frac{n}2\big|<x{\frac{\sqrt n}2}\Big)
           \,=\, 
   \int_{-x}^{ x }
   e^{- \tfrac{s^2}{2} }
   %\exp\Big(- \frac{s^2}{2} \Big)\,
   \frac{ds}{
             \sqrt{{2}{\pi}}
   }
  \,.
\end{equation}
%It is then an easy matter to deduce the 
%de Moivre theorem 
It is then an easy matter to obtain the result
for an arbitrary interval: just apply the previous result to two intervals 
$[-x,x]$ and $[-y,y]$, with $x<y$, 
subtract the two results and use the symmetry of 
the binomial and Gaussian distributions.
\subsection{Avoiding Stirling's formula.}
Let us discuss now the role played by Stirling's formula.
In the previous argument, it is used only to prove the limit~\eqref{fstir}.
However this limit can be proved directly with the help of 
the Wallis integrals. These are defined as 
\begin{equation}
  \forall n\geq 0\qquad 
  W_n\,=\,\int_0^{\pi/2}(\sin t )^n\,dt\,.  
\end{equation}
Let us recall some classical facts on the sequence $(W_n)_{n\geq 0}$.
Since the sinus takes its values 
in $[0,1]$ on the interval $[0,\pi/2]$, 
then
\begin{equation}
\label{wrec}
  \forall n\geq 1\qquad 
  %W_{n+1}\,\leq\,W_n\,.
  W_{n+1}\,\leq\,W_n\,\leq\, W_{n-1}\,.
\end{equation}
An integration by parts yields the classical recurrence formula
\begin{equation}
  \label{wipp}
  \forall n\geq 2\qquad 
  nW_n\,=\, (n-1)W_{n-2}\,.
\end{equation}
Multiplying~\eqref{wipp} by $W_{n-1}$, 
we obtained the remarkable fact that 
\begin{equation}
  \label{reame}
  \forall n\geq 2\qquad 
  nW_nW_{n-1}\,=\, W_1W_0\,=\,\frac{\pi}{2}\,.
\end{equation}
We can compute the exact values of the integrals
by iterating formula~\eqref{wipp}:
%and we obtain
\begin{gather}
  \label{wdn1}
  W_{2n}\,=\,
  \frac{(2n-1)\times\cdots\times 1}{(2n)\times\cdots\times 2}\frac{\pi}{2}\,,\\
  \label{wdn2}
W_{2n+1}\,=\,
  \frac{(2n)\times\cdots\times 2}
  {(2n+1)\times\cdots\times 3}\,.
\end{gather}
It follows also from~\eqref{wipp} that
$W_{n+1}\sim W_{n-1}$ as $n$ goes to $\infty$. This fact, 
combined with the inequalities~\eqref{wrec}, yields that
$W_{n}\sim W_{n+1}$ as $n$ goes to $\infty$. 
It follows then from~\eqref{reame} that
$n(W_{n})^2\sim\pi/2$, whence 
\begin{equation}
  \label{equw}
  W_n\,\sim\, \sqrt{\frac{\pi}{2n}}
  \quad\text{as}\quad  n\to\infty \,.
\end{equation}
The Wallis integrals are usually employed to prove 
Stirling's formula.
%the Stirling formula.
In a first step,
one proves that $n!\sim C(n/e)^n\sqrt{n}$, where $C$ 
is a positive constant. In a second step, one inserts 
this estimate in the formula for $W_{2n}$ and uses the 
estimate~\eqref{equw} to conclude that $C=\sqrt{2\pi}$.
However, if one only wishes to prove the 
limit~\eqref{fstir},
then these two steps,
and in fact all the asymptotic computations of the equivalents,
can be completely bypassed.
%If one wishes only to prove the limit~\eqref{fstir},
%then we can completely bypass these two steps, in fact even 
%all the asymptotic computations of the equivalents.
The trick consists in rewriting the formulas~\eqref{wdn1} and~\eqref{wdn2}
as product of factorials:
\begin{comment}
  \label{rwdn1}
  W_{2n}\,=\,
  \frac{(2n)!} {2^{2n}(n!)^2}
  \frac{\pi}{2}\,,\\
  \label{rwdn2}
W_{2n+1}\,=\,
\frac{1}{2n+1}
  \frac{2^{2n}(n!)^2}{(2n)!} \,.
\end{comment}
\begin{equation}
  \label{rwdn}
  W_{2n}\,=\,
  \frac{(2n)!} {2^{2n}(n!)^2}
  \frac{\pi}{2}\,,\qquad
W_{2n+1}\,=\,
\frac{1}{2n+1}
  \frac{2^{2n}(n!)^2}{(2n)!} \,.
\end{equation}
The central probability $P(S_{2n}=n)$ pops up naturally and we see that
\begin{equation}
  \label{pwdn}
  W_{2n}\,=\,
 % \frac{(2n)!} {2^{2n}(n!)^2}
  P(S_{2n}=n)
  \frac{\pi}{2}
  \,,\qquad
W_{2n+1}\,=\,
\frac{1}{(2n+1)
  P(S_{2n}=n)
}
  \,.
\end{equation}
\begin{comment}
  \label{spwdn}
  W_{2n}\,=\,
 % \frac{(2n)!} {2^{2n}(n!)^2}
  P(S_{2n}=n)
  \frac{\pi}{2}
  \,,\qquad
W_{2n+1}\,=\,
\frac{1}{(2n+2)
  P(S_{2n+1}=n)
}
  \,.
\end{comment}
From there, we could easily deduce the limit~\eqref{fstir}
from the asymptotic equivalent~\eqref{equw}.
%However, we can even do simpler and better. Indeed, 
However, we can do it even more simply and efficiently. Indeed, 
let us look 
at the ratio 
$W_{2n}/ W_{2n+1}$.
On one hand,
from the expressions~\eqref{pwdn}, we have
\begin{equation}
  \label{lwdn}
  \frac{ W_{2n}}{ W_{2n+1}}\,=\,
  \Big(n+\frac12\Big)\pi\big(P(S_{2n}=n)\big)^2 \,.
  %(2n+1)\big(P(S_{2n}=n)\big)^2 \frac{\pi}{2}\,.
\end{equation}%with the help of~\eqref{equw}.
On the other hand, it follows from the inequalities~\eqref{wrec}
and the recurrence relation~\eqref{wipp} that
\begin{equation}
  \label{zwdn}
  1\,\leq\,
  \frac{ W_{2n}}{ W_{2n+1}}\,\leq\,
  \frac{ W_{2n-1}}{ W_{2n+1}}\,=\,
  \frac{ {2n+1}}{ {2n}}\,=\,
  1+\frac{ {1}}{ {2n}}
  \,.
\end{equation}
Combining~\eqref{lwdn} and~\eqref{zwdn}, we conclude that,
for any $n\geq 1$, 
\begin{equation}
  \label{cwdn}
  %\forall n\geq 1\qquad 
  1\,\leq\,
  \sqrt{\Big(n+\frac12\Big)\pi}P(S_{2n}=n)
  %(2n+1)\big(P(S_{2n}=n)\big)^2 \frac{\pi}{2}\,.
  \,\leq\,
  \sqrt{1+\frac{ {1}}{ {2n}}}
  \,,
\end{equation}
which we rewrite in the more convenient form 
\begin{equation}
  \label{mwdn}
  %\forall n\geq 1\qquad 
  \frac{1}{\sqrt{1+\frac1{2n}}}
  \,\leq\,
  \sqrt{n\pi}
  P(S_{2n}=n)
  %(2n+1)\big(P(S_{2n}=n)\big)^2 \frac{\pi}{2}\,.
  \,\leq\,1
  %\sqrt{1+\frac{ {1}}{ {2n}}}
  \,.
\end{equation}
The inequalities~\eqref{mwdn} readily yield the desired limit~\eqref{fstir}.
%as well as its counterpart 
For the odd integers, 
we notice that
\begin{equation}
  P(S_{2n+1}=n)\,=\,
  \frac{2n+1}{2n+2}
  P(S_{2n}=n)\,,
\end{equation}
and we deduce from~\eqref{mwdn} that
\begin{equation}
  \label{dwdn}
  %\forall n\geq 1\qquad 
  %\frac{1}{\sqrt{1+\frac1{2n}}}
  \frac{1}{\Big({1+\frac1{2n}}\Big)^{3/2}}
  \,\leq\,
  \sqrt{n\pi}
  P(S_{2n+1}=n)
  %(2n+1)\big(P(S_{2n}=n)\big)^2 \frac{\pi}{2}\,.
  \,\leq\,1
  %\sqrt{1+\frac{ {1}}{ {2n}}}
  \,.
\end{equation}
Sending $n$ to $\infty$ in~\eqref{dwdn}, we obtain
%These inequalities yield 
%as well as its 
the counterpart 
of the limit~\eqref{fstir}
for the odd integers. 
\subsection{Non-asymptotic inequalities.}
One merit of our proof is that it produces interesting
non-asymptotic inequalities. 
We restate next these inequalities, and we combine them with the
inequalities coming from the Wallis integrals.
%
%that are obtained in the course 
%We shall finally restate
%the non-asymptotic inequalities that are obtained in the course 
%of our proof.
%We shall next combine lemma~\ref{nacen} with the 
%inequalities we derived to prove de Moivre's theorem.
%For the upper bound
Let $n\geq 1$ and $x>0$ be such that $n\geq x^2$.
Inequalities~\eqref{uwex} and~\eqref{mwdn} yield
\begin{align}
  \label{aakwex}
  %\forall n\geq 1\quad \forall n\geq x^2\qquad\\
P\Big(\big|S_{2n}-{n}\big|<x\sqrt{\frac{n}{2}}\Big)
  %P_{2n}(x)
           &\,\leq\, 
           %\sum_{\substack{k\in\zun\\ |k-n|< x\sqrt{\frac{n}2}}}
           %\sqrt{2n}P(S_{2n}=n)
  %\Big({1+\frac{ {1}}{ {4n}}}\Big)
             %e^{ \frac{x^3}{2\sqrt{2n}}+\frac{ {1}}{ {4n}}}
           e^{ \frac{x^3}{2\sqrt{2n}}}
   \int_{-x}^{ x }
   e^{- \tfrac{s^2}{2} }
   %\exp\Big(- \frac{s^2}{2} \Big)\,
   \frac{ds}{
             \sqrt{{2}{\pi}}
   }
   %\times \exp\Big( -\frac{x}{\sqrt{2n}}
             +\frac{1}{ \sqrt{{n\pi}} }\cr
           &\,\leq\, 
           %\sum_{\substack{k\in\zun\\ |k-n|< x\sqrt{\frac{n}2}}}
           %\sqrt{2n}P(S_{2n}=n)
  %\Big({1+\frac{ {1}}{ {4n}}}\Big)
             %e^{ \frac{x^3}{2\sqrt{2n}}+\frac{ {1}}{ {4n}}}
   \int_{-x}^{ x }
   e^{- \tfrac{s^2}{2} }
   %\exp\Big(- \frac{s^2}{2} \Big)\,
   \frac{ds}{
             \sqrt{{2}{\pi}}
   }
   %\times \exp\Big( -\frac{x}{\sqrt{2n}}
           +e^{ \frac{x^3}{2\sqrt{2n}}}-1
             +\frac{1}{ \sqrt{{n\pi}} }
  \,,
\end{align}
because the integral is less than $1$.
%Using the inequality~\eqref{bce}
Using inequality~\eqref{bce}
with $t=1/\sqrt{n\pi}$,
we have
\begin{equation}
  \label{iste}
           e^{ \frac{x^3}{2\sqrt{2n}}}
             +\frac{1}{ \sqrt{{n\pi}} }
             \,\leq\,
           e^{ \frac{x^3}{2\sqrt{2n}}}
             \Big(1+\frac{1}{ \sqrt{{n\pi}} }\Big)
             \,\leq\,
           e^{ \frac{x^3}{2\sqrt{2n}}
             +\frac{1}{ \sqrt{{n\pi}} }
           }
             \,\leq\,
           e^{ \frac{x^3+2}{2\sqrt{2n}}}\,.
\end{equation}
%we obtain that
%we rewrite it in the following more convenient form:
Combining the inequalities~\eqref{aakwex} and~\eqref{iste}, 
we obtain
\begin{multline}
  \label{nakwex}
  \forall x>0\quad
\forall n\geq \max(x^2,1)\qquad\cr %\\
P\Big(\big|S_{2n}-{n}\big|<x\sqrt{\frac{n}{2}}\Big)
  %P_{2n}(x)
           \,\leq\, 
           %\sum_{\substack{k\in\zun\\ |k-n|< x\sqrt{\frac{n}2}}}
           %\sqrt{2n}P(S_{2n}=n)
  %\Big({1+\frac{ {1}}{ {4n}}}\Big)
             %e^{ \frac{x^3}{2\sqrt{2n}}+\frac{ {1}}{ {4n}}}
   \int_{-x}^{ x }
   e^{- \tfrac{s^2}{2} }
   %\exp\Big(- \frac{s^2}{2} \Big)\,
   \frac{ds}{
             \sqrt{{2}{\pi}}
   }
   %\times \exp\Big( -\frac{x}{\sqrt{2n}}
           +e^{ \frac{x^3+2}{2\sqrt{2n}}}-1
             %+\frac{1}{ \sqrt{{n}} }
  \,.
\end{multline}
We try now to derive a similar lower bound.
Let $n\geq 1$ and $x>0$ be such that $n\geq 2x^2$.
Inequalities~\eqref{fvex} and~\eqref{mwdn} yield
\begin{align}
  \label{iuvex}
  %\forall n\geq 1\quad \forall n\geq 2x^2\qquad\\
%  P_{2n}(x)
P\Big(\big|S_{2n}-{n}\big|<x\sqrt{\frac{n}{2}}\Big)
           &\,\geq\, 
           %\sum_{\substack{k\in\zun\\ |k-n|< x\sqrt{\frac{n}2}}}
           %\sqrt{2n}P(S_{2n}=n)
  \frac{1}{\sqrt{1+\frac1{2n}}}
  %\Big({1-\frac{ {1}}{ {4n}}}\Big)
           \Bigg(
             e^{ -\frac{x^3}{\sqrt{2n}}}
   \int_{-x}^{ x }
   %\exp\Big(- \frac{s^2}{2} \Big)\,
   e^{- \tfrac{s^2}{2} }
   \frac{ds}{
             \sqrt{{2}{\pi}}
   }
   %\times \exp\Big( -\frac{x}{\sqrt{2n}}
             -
             \frac{1}{
             \sqrt{{n}{\pi}}
             }
           \Bigg)
           \cr
           &\,\geq\, 
           %\sum_{\substack{k\in\zun\\ |k-n|< x\sqrt{\frac{n}2}}}
           %\sqrt{2n}P(S_{2n}=n)
  \frac{1}{\sqrt{1+\frac1{2n}}}
  %\Big({1-\frac{ {1}}{ {4n}}}\Big)
             e^{ -\frac{x^3}{\sqrt{2n}}}
   \int_{-x}^{ x }
   %\exp\Big(- \frac{s^2}{2} \Big)\,
   e^{- \tfrac{s^2}{2} }
   \frac{ds}{
             \sqrt{{2}{\pi}}
   }
   %\times \exp\Big( -\frac{x}{\sqrt{2n}}
             -
             \frac{1}{
             \sqrt{{2n}}
             }
           \,.
\end{align}
Applying the inequality~\eqref{bce}
with $t=1/(2n)$,
we have 
\begin{equation}
  \label{fhh}
  \frac{1}{\sqrt{1+\frac1{2n}}}
             e^{ -\frac{x^3}{\sqrt{2n}}}
             \,\geq\, 
             e^{ 
  -\frac{ {1}}{ {4n}}
               -\frac{x^3}{\sqrt{2n}}
}\,.
             %\,\geq\,
  %{1-\frac{ {1}}{ {4n}}}
             %{ -\frac{x^3}{\sqrt{2n}}}\,.
\end{equation}
Substituting~\eqref{fhh} into~\eqref{iuvex},
and using the fact that
the integral is less than $1$,
%Combining the two previous inequalities, 
we obtain 
\begin{equation}
  \label{juvex}
  %\forall n\geq 1\quad \forall n\geq 2x^2\qquad\\
%  P_{2n}(x)
P\Big(\big|S_{2n}-{n}\big|<x\sqrt{\frac{n}{2}}\Big)
           \,\geq\, 
           %\sum_{\substack{k\in\zun\\ |k-n|< x\sqrt{\frac{n}2}}}
           %\sqrt{2n}P(S_{2n}=n)
   \int_{-x}^{ x }
   %\exp\Big(- \frac{s^2}{2} \Big)\,
   e^{- \tfrac{s^2}{2} }
   \frac{ds}{
             \sqrt{{2}{\pi}}
   }
   %\times \exp\Big( -\frac{x}{\sqrt{2n}}
             +e^{ 
  -\frac{ {1}}{ {4n}}
               -\frac{x^3}{\sqrt{2n}}
             }
             -1
             -
             \frac{1}{
             \sqrt{{2n}}}
           \,.
  \end{equation}
  Moreover, we have,
thanks to inequality~\eqref{bce} applied
with $t=1/\sqrt{2n}$,
  \begin{align}
    \label{kuw}
             e^{ -\frac{ {1}}{ {4n}} -\frac{x^3}{\sqrt{2n}} }
             -1-
             \frac{1}{ \sqrt{{2n}}}&\,\geq\,
             e^{ -\frac{ {1}}{ {4n}} -\frac{x^3}{\sqrt{2n}} }
             -
             e^{\frac{1}{ \sqrt{{2n}}}}\cr
             %\,\geq\,
             %1-e^{ \frac{ x^3+2}{\sqrt{n}} +
             %\frac{1}{ \sqrt{{n}{\pi}}}}
                &\,\geq\,
             1-e^{ \frac{ {1}}{ {4n}}+ \frac{x^3}{\sqrt{2n}} +
             \frac{1}{ \sqrt{{2n}}}}
             \,\geq\,
             1-e^{ \frac{ x^3+2}{\sqrt{n}}}\,.
  \end{align}
  We conclude from~\eqref{juvex} and~\eqref{kuw} that
%we obtain after some minor modifications
\begin{multline}
  \label{uuvex}
  \forall x>0\quad
\forall n\geq \max(2x^2,1)\qquad\cr %\\
  %\forall n\geq 1\quad \forall n\geq 2x^2\quad%\\
%  P_{2n}(x)
P\Big(\big|S_{2n}-{n}\big|<x\sqrt{\frac{n}{2}}\Big)
           \,\geq\, 
           %\sum_{\substack{k\in\zun\\ |k-n|< x\sqrt{\frac{n}2}}}
           %\sqrt{2n}P(S_{2n}=n)
  %\Big({1-\frac{ {1}}{ {4n}}} { -\frac{x^3}{\sqrt{2n}}} \Big)
   \int_{-x}^{ x }
   e^{- \tfrac{s^2}{2} }
   %\exp\Big(- \frac{s^2}{2} \Big)\,
   \frac{ds}{
             \sqrt{{2}{\pi}}
   }
   %\times \exp\Big( -\frac{x}{\sqrt{2n}}
             +1-e^{ \frac{ x^3+2}{\sqrt{n}}}\,.
  %{ -\frac{x^3}{\sqrt{2n}}} 
\end{multline}
We can finally combine
the inequalities~\eqref{nakwex} 
and~\eqref{uuvex} into a single inequality: 
%to get the following inequality:
\begin{multline}
  \label{luvex}
  %\forall n\ge\quad \forall n\geq 2x^2\quad%\\
  \forall x>0\quad
\forall n\geq \max(2x^2,1)\qquad\cr
%  P_{2n}(x)
\Big|
P\Big(\big|S_{2n}-{n}\big|<x\sqrt{\frac{n}{2}}\Big)
-
   \int_{-x}^{ x }
   e^{- \tfrac{s^2}{2} }
   %\exp\Big(- \frac{s^2}{2} \Big)\,
   \frac{ds}{
             \sqrt{{2}{\pi}}
   }
\Big|\,\leq\,
   %\times \exp\Big( -\frac{x}{\sqrt{2n}}
           %e^{ \tfrac{x^3+2}{2\sqrt{n}}}-1
           e^{ \tfrac{x^3+2}{\sqrt{n}}}-1
  %{ -\frac{x^3}{\sqrt{2n}}} 
           \,.
\end{multline}%The literature is full of more sophisticated 
For the odd integers, 
with the help of the inequalities~\eqref{oouwvex}
and~\eqref{dwdn}, 
we obtain the following similar inequality (the details are omitted):
\begin{multline}
  \label{siluvex}
  %\forall n\ge\quad \forall n\geq 2x^2\quad%\\
  \forall x>0\quad
\forall n\geq \max(2x^2,1)\qquad%\\
  \cr
%  P_{2n}(x)
 % P_{2n+1}(x)\,=\, 
\Big|
  P\Big(\big|S_{2n+1}-\frac{2n+1}{2}\big|<
  \frac{x}{2}\sqrt{2n+1}
\Big)
-
   \int_{-x}^{ x }
   e^{- \tfrac{s^2}{2} }
   %\exp\Big(- \frac{s^2}{2} \Big)\,
   \frac{ds}{
             \sqrt{{2}{\pi}}
   }
\Big|\,\leq\,
   %\times \exp\Big( -\frac{x}{\sqrt{2n}}
           %e^{ \tfrac{x^3+2}{2\sqrt{n}}}-1
           e^{ \tfrac{x^3+2}{\sqrt{n}}}-1
  %{ -\frac{x^3}{\sqrt{2n}}} 
           \,.
\end{multline}
Combining the inequalities~\eqref{luvex} and~\eqref{siluvex},
we finally obtain the inequality~\eqref{gluvex} stated in the 
introduction.
\begin{comment}
  \label{oluvex}
  \forall n\geq 0\quad
\forall n\geq 2x^2\qquad\\
%  P_{2n}(x)
\Big|
P\Big(\big|S_{2n+1}-\frac{2n+1}2\big|<\frac{x}{2}\sqrt{{2n+1}}\Big)
-
   \int_{-x}^{ x }
   e^{- \tfrac{s^2}{2} }
   %\exp\Big(- \frac{s^2}{2} \Big)\,
   \frac{ds}{
             \sqrt{{2}{\pi}}
   }
\Big|\,\leq\,
   %\times \exp\Big( -\frac{x}{\sqrt{2n}}
           e^{ \tfrac{x^3+2}{2\sqrt{n}}}-1
  %{ -\frac{x^3}{\sqrt{2n}}} 
           \,.
\end{comment}
\section{de Moivre--Laplace theorem.}
\label{dmlt}
We consider 
%now 
the general case $p\in ]0,1[$.
The nice symmetries of the case $p=1/2$ disappear.
To start the computations, we follow the path laid out 
by Feller 
(see \cite{F}, section VII.3, pages 182 and~183). 
First we set $q=1-p$, and we introduce the index $m$ of the central
term of the binomial distribution $B(n,p)$ as the unique 
integer $m$ such that 
\begin{equation}
 m=np+\delta\,,\qquad -q<\delta\leq p\,. 
\end{equation}
Let $x>0$ and let us define 
%. We shall estimate the probability 
\begin{equation}
  \forall n\geq 1\qquad 
  P_n(x)\,=\, P\big(m\leq S_n<m+x{{\sqrt {npq}}}\big)\,.
\end{equation}
We fix $n\geq 1$ and we start by 
writing down the expression of $P_{n}(x)$:
\begin{equation}
  \label{gdex}
  P_{n}(x)
           \,=\, 
           %\sum_{\substack{k\in\zun\\ |k-n|< x\sqrt{\frac{n}2}}}
           \kern-10pt
           %\sum_{\substack{k\in\zun\\ m\leq k<m+ x\sqrt{{npq}}}}
           \sum_{{m\leq k<m+ x\sqrt{{npq}}}}
           \kern-15pt
           P\big(S_{n}=k\big)
           \,=\, 
           \kern-10pt
           %\sum_{\substack{k\in\zun\\ |k-n|< x\sqrt{\frac{n}2}}}
           %\sum_{\substack{k\in\zun\\ m\leq k<m+ x\sqrt{{npq}}}}
           \sum_{{m\leq k<m+ x\sqrt{{npq}}}}
           \binom{n}{k}p^kq^{n-k}\,.
\end{equation}
We make the change of variable $j=k-m$, so that~\eqref{gdex} becomes
\begin{equation}
  \label{grex}
  P_{n}(x)
           \,=\, 
           %\sum_{\substack{k\in\zun\\ |k-n|< x\sqrt{\frac{n}2}}}
           %\sum_{\substack{k\in\zun\\ 0\leq j< x\sqrt{{npq}}}}
           \sum_{{0\leq j< x\sqrt{{npq}}}}
           \binom{n}{m+j}p^{m+j}q^{n-m-j}\,.
\end{equation}
We proceed as in the symmetric case. The first step consists 
in rewriting the binomial coefficient 
           $\binom{n}{m+j}$ as the product of 
           $\binom{n}{m}$ and a fraction.
For $j\geq 1$, we have 
\begin{align}
  \label{gact1}
  (m+j)!&\,=\,m!\,(m+1)\cdots(m+j)\,,\\ \label{gact2}
  (n-m-j)!&\,=\,\frac{(n-m)!}{(n-m)\cdots (n-m-j+1)}\,,
\end{align}
whence
\begin{equation}
  \label{dent}
           \binom{n}{m+j}\,=\,
           \frac{n!}{{(m+j)!}{(n-m-j)!}}
           \,=\,
           %\frac{(2n)!}{n! n!}
           \binom{n}{m}
%\frac{(n-m-j+1)\cdots (n-m)}
\frac{(n-m)\cdots (n-m-j+1)}
  {(m+1)\cdots(m+j)}
  \,.
\end{equation}
We introduce next a notation for the product appearing in~\eqref{dent},
and we incorporate as well the factor $(p/q)^j$ coming from~\eqref{grex} :
\begin{equation}
  \label{qent}
  \forall j\geq 1\qquad
  \pi(j)
           \,=\,
\frac{(n-m)\cdots (n-m-j+1)}
%\frac{(n-m-j+1)\cdots (n-m)}
  {(m+1)\cdots(m+j)}
  \Big(
  \frac{p}{q}
  \Big)^j
           \,=\,
           \prod_{1\leq\ell\leq j}
  \frac{(n-m-\ell+1)p}{(m+\ell)q}
  %\frac{n-m+\ell-j}{m+\ell}
  \,.
\end{equation}
\begin{comment}
  \label{qent}
  \forall j\geq 1\qquad
  \pi(j)
           &\,=\,
\frac{(n-m)\cdots (n-m-j+1)}
%\frac{(n-m-j+1)\cdots (n-m)}
  {(m+1)\cdots(m+j)}
  \Big(
  \frac{p}{q}
  \Big)^j\cr
           &\,=\,
           \prod_{1\leq\ell\leq j}
  \frac{(n-m-\ell+1)p}{(m+\ell)q}
  %\frac{n-m+\ell-j}{m+\ell}
  \,.
\end{comment}
We set also $\pi(0)=1$.
With this notation, the formula~\eqref{grex} can be rewritten as
\begin{equation}
  \label{bwex}
  P_{n}(x)
           \,=\, 
           %\sum_{\substack{k\in\zun\\ |k-n|< x\sqrt{\frac{n}2}}}
           %\frac{2}{2^{2n}} \binom{2n}{n}
           P(S_{n}=m)
           \sum_{{0\leq j< x\sqrt{{npq}}}}
          % \binom{n}{m+j}
           %\pi(j)p^{j}q^{-j}\,.
           \pi(j)\,.
           %\Big(\frac{p}{q}\Big)^{j}\,.
\end{equation}
We focus next on the product $\pi(j)$. Replacing $m$ by $np+\delta$,
we obtain the following expression:
\begin{equation}
  \label{qrnt}
  \forall j\geq 1\qquad
  \pi(j)
           %\,=\,
           %\prod_{1\leq\ell\leq j}
  %\frac{nq-\delta-\ell+1}{np+\delta+\ell}
           \,=\,
           \prod_{1\leq\ell\leq j}
  \frac{npq-(\delta+\ell-1)p}{npq+(\delta+\ell)q}
  \,.
\end{equation}%\subsection{lower bound of $\pi(j,n)$.}
%\subsection{Upper bound of $\pi(j,n)$.}
\subsection{Upper bound.}
Let us fix $j\in\zun$. To prepare for the upper bound, we rewrite 
$\pi(j)$ as follows:
\begin{equation}
  \label{vanop}
  \pi(j) 
  %\,=\, \prod_{1\leq\ell\leq j}
  %\frac{npq-(\delta+\ell-1)p}{npq+(\delta+\ell)q}
   \,=\, \prod_{1\leq\ell\leq j}
   \Big(1-
  \frac{\delta+\ell-p}{npq+(\delta+\ell)q}
     \Big)
  \,,
\end{equation}
and we use inequality~\eqref{bce} 
to deduce from~\eqref{vanop} that
\begin{multline}
  \label{wanor}
  \pi(j) \,\leq\, 
  %\prod_{1\leq\ell\leq j}
  %\exp\Big(-
  %\frac{j}{n+\ell}
    %\Big)
   %\,=\, 
  \exp\Big(
  - \sum_{1\leq\ell\leq j}
  \frac{\delta+\ell-p}{npq+(\delta+\ell)q}
     \Big)\cr
   \,\leq\, 
  \exp\Big(
  - \sum_{1\leq\ell\leq j}
  \frac{\ell-1}{npq+pq+jq}
     \Big)
     %\cr
   \,=\, 
  \exp\Big(
  - 
  \frac{j(j-1)}{2\big((n+1)pq+jq\big)}
     \Big)
  \,.
\end{multline}
Let us remark that the previous computations are more complicated 
than the ones we did in the symmetric case, not only because of the
presence of the parameters $p,q$, but because there is a specific
cancellation when $p=q$.
Indeed, suppose that we try to proceed exactly as in the symmetric case.
Instead of 
formulas~\eqref{qrnt}
and~\eqref{vanop}, we would start from the 
equivalent formula where $\ell$
is changed into $j-\ell+1$ in the numerator, and we would get:
\begin{equation}
  \label{yanop}
  \pi(j) \,=\, 
           \prod_{1\leq\ell\leq j}
  \frac{npq-(\delta+j-\ell)p}{npq+(\delta+\ell)q}
   \,=\, \prod_{1\leq\ell\leq j}
   \Big(1-
  \frac{\delta+\ell(q-p)+jp}{npq+(\delta+\ell)q}
     \Big)
  \,.
\end{equation}
In the symmetric case, the term $\ell (q-p)$ in the last expression vanished,
and this leads to a simpler upper bound. In the non-symmetric case, one has 
to go through an intermediate 
exponential upper bound and perform the summation over $\ell$ 
in the argument of the exponential, as we did in~\eqref{wanor}.
Moreover, these  steps are more transparent when starting from 
formula~\eqref{vanop} than with formula~\eqref{yanop}, and to 
tell the truth, 
the product in the numerator of 
formula~\eqref{qrnt} was 
already written in reverse order 
compared to~\eqref{pent}.

Coming back to~\eqref{wanor},
we use inequality~\eqref{ine1}
to obtain that, for $0\leq j\leq (n+1)p$,
\begin{multline}
  \label{gbnos}
  \frac{j(j-1)}{2\big((n+1)pq+jq\big)}
   \,=\,
\frac{j(j-1)}{2(n+1)pq}
   %\frac{1}{ \displaystyle 1+\frac{j}{(n+1)p}
   \frac{1}{ 1+\frac{j}{(n+1)p}
   }\cr
   \,\geq\,
\frac{j(j-1)}{2(n+1)pq}
   \Big( 1-\frac{j}{(n+1)p}\Big)
   \,\geq\,
\frac{(j-1)^2}{2(n+1)pq}
   -\frac{j^3}{2n^2p^2q}
  \,.
\end{multline}
Substituting~\eqref{gbnos} in~\eqref{wanor},
%and~\eqref{vanop}, we conclude that
we conclude that
\begin{equation}
  \label{vinos}
\forall j\,,\quad 0\leq j\leq (n+1)p\,,\qquad
  \pi(j) 
  \,\leq \, 
   \exp\Big(-
\frac{(j-1)^2}{2(n+1)pq}
   +\frac{j^3}{2n^2p^2q}
   \Big)
  \,.
\end{equation}
\begin{comment}
Contrary to the previous case, this inequality is as good as the one 
obtained by 
Feller 
(see \cite{F}, page 180, between formulas~$(2.3)$ and $(2.4)$).
\end{comment}
%\subsection{Upper bound of $P_{2n}(x)$.}
With the upper bound of $\pi(j)$ in hand, we will compute 
an upper bound on 
$P_{n}(x)$. Let us fix $x>0$ and let $n$ be an integer 
such that $n\geq (q/p)x^2$. This ensures that 
\begin{equation}
  \label{lfnst}
           {{0\leq j< x\sqrt{{npq}}}}
           \quad
           \Longrightarrow
           \quad
           %{j}\leq (n+1)p\,,\quad
           {j}\leq np\,,\quad
         %  \frac{j}{n}\leq\frac12\,,\quad
   \frac{j^3}{2n^2p^2q}
   %\leq\frac{x^3}{2}\sqrt{\frac{q}{pn}}\,.
   \leq\frac{x^3}{\sqrt{{pn}}}\,.
\end{equation}
%In particular, for these values of $j$, the inequality~\eqref{vinos}
In particular, for these values of $j$, the inequality~\eqref{vinos}
holds as well, whence
%and we have, thanks to~\eqref{lfnst},
\begin{equation}
  \label{bvinos}
           {{0\leq j< x\sqrt{{npq}}}}
           \quad
           \Longrightarrow
           \quad 
  \pi(j) 
  \,\leq \, 
   \exp\Big(-
\frac{(j-1)^2}{2(n+1)pq}
   %+\frac{x^3}{2}\sqrt{\frac{q}{pn}}
   +\frac{x^3}{\sqrt{{pn}}}
   \Big)
  \,.
\end{equation}
%we can substitute inequality~\eqref{linos}
%In particular, we can substitute inequality~\eqref{linos}
Substituting inequality~\eqref{bvinos}
into the sum appearing in~\eqref{bwex}, we get
\begin{equation}
  \label{twex}
  P_{n}(x)
           \,\leq\, 
           %\sum_{\substack{k\in\zun\\ |k-n|< x\sqrt{\frac{n}2}}}
           %\frac{2}{2^{2n}} \binom{2n}{n}
           P(S_{n}=m)
           \sum_{{0\leq j< x\sqrt{{npq}}}}
          % \binom{n}{m+j}
           %\pi(j)p^{j}q^{-j}\,.
%           \Big(\frac{p}{q}\Big)^{j}
   \exp\Big(-
\frac{(j-1)^2}{2(n+1)pq}
   %+\frac{x^3}{2}\sqrt{\frac{q}{pn}}
   +\frac{x^3}{\sqrt{{pn}}}
   \Big)
           \,.
\end{equation}
%So, why did we bother to rederive 
%inequality~\eqref{tmyin} from scratch? 
%It should be noted that the proof of Keller and Kindler
%is quite complex and difficult (although not as mysterious
%as Talagrand's proof), and it involves also other
%deep results, typically hypercontractivity estimates.
Since the function $\exp(-t)$ is decreasing on $[0,+\infty[$, then
\begin{equation}
  \label{amdosi}
  \forall j\geq 2\qquad 
   \exp\Big(- 
\frac{(j-1)^2}{2(n+1)pq}
   \Big)\,\leq\,
   \int_{j-2}^{j-1}
   e^{- 
\frac{t^2}{2(n+1)pq}
   }
  % \exp\Big(- \frac{t^2}{n} \Big)\,
   dt\,.
\end{equation}
We sum the inequality~\eqref{amdosi} over $j$
           such that $0\leq j< x\sqrt{{npq}}$, and we bound 
the first two terms by $1$:
\begin{equation}
  \label{samdosi}
           \sum_{{0\leq j< x\sqrt{{npq}}}}
   \exp\Big(- 
\frac{(j-1)^2}{2(n+1)pq}
   \Big)\,\leq\,
  2+
  %\int_0^{ \lfloor x\sqrt{{npq}}\rfloor-1}
  \int_0^{  x\sqrt{{npq}}}
  e^{- \frac{t^2}{2(n+1)pq} }
   %\exp\Big(- \frac{t^2}{2(n+1)pq} \Big)
   \,dt\,.
\end{equation}
Substituting inequality~\eqref{samdosi}
%into the sum appearing in~\eqref{twex}, we get
into~\eqref{twex}, we get
\begin{equation}
  \label{awex}
  P_{n}(x)
           \,\leq\, 
           %\sum_{\substack{k\in\zun\\ |k-n|< x\sqrt{\frac{n}2}}}
           %\frac{2}{2^{2n}} \binom{2n}{n}
           P(S_{n}=m)
   %\exp\Big( \frac{x^3}{2}\sqrt{\frac{q}{pn}}\Big)
 %e^{ \frac{x^3}{2}\sqrt{\frac{q}{pn}}}
 \,
 e^{ \frac{x^3}{\sqrt{{pn}}} }
          % \binom{n}{m+j}
           %\pi(j)p^{j}q^{-j}\,.
%           \Big(\frac{p}{q}\Big)^{j}
\Bigg(2+
  \int_0^{ x\sqrt{{npq}}}
  e^{- \frac{t^2}{2(n+1)pq} }
   %\exp\Big(- \frac{t^2}{2(n+1)pq} \Big)
   \,dt
   \Bigg)
           \,.
\end{equation}
We make the change of variable 
$s { \sqrt{{(n+1)pq}} }=t$ in~\eqref{awex} and the inequality becomes
\begin{equation}
  \label{aoex}
  P_{n}(x)
           \,\leq\, 
           %\sum_{\substack{k\in\zun\\ |k-n|< x\sqrt{\frac{n}2}}}
           %\frac{2}{2^{2n}} \binom{2n}{n}
           P(S_{n}=m)
   %\exp\Big( \frac{x^3}{2}\sqrt{\frac{q}{pn}}\Big)
 %e^{ \frac{x^3}{2}\sqrt{\frac{q}{pn}}}
 \,
 e^{ \frac{x^3}{\sqrt{{pn}}} }
          % \binom{n}{m+j}
           %\pi(j)p^{j}q^{-j}\,.
%           \Big(\frac{p}{q}\Big)^{j}
\Bigg(2+
  %\int_0^{ x\sqrt{{npq}}}
 \sqrt{{(n+1)pq} }
  \int_0^{ x}
  e^{- \frac{s^2}{2} }
   %\exp\Big(- \frac{t^2}{2(n+1)pq} \Big)
   \,ds
   \Bigg)
           \,.
\end{equation}
We finally factorize 
$ \sqrt{{npq} }$,
we 
introduce the normalizing factor $\sqrt{2\pi}$
and we use the inequalities
%${ \sqrt{{1+\frac1n} }\leq e^{1/n}}$ to get
\begin{equation}
  \sqrt{{1+\frac1n} }\,\leq \,
  1+\frac1{2n} \,\leq \,{e^{\frac1n}}
  %\,\leq \,{e^{\frac1{2n}}}
  \,\leq \,
 e^{ \frac{1}{\sqrt{{pn}}} }
\end{equation}
to get
\begin{equation}
  \label{azoex}
  P_{n}(x)
           \,\leq\, 
           %\sum_{\substack{k\in\zun\\ |k-n|< x\sqrt{\frac{n}2}}}
           %\frac{2}{2^{2n}} \binom{2n}{n}
 \sqrt{{2\pi npq} }\,
           P(S_{n}=m)
   %\exp\Big( \frac{x^3}{2}\sqrt{\frac{q}{pn}}\Big)
 %e^{ \frac{x^3}{2}\sqrt{\frac{q}{pn}}}
 \,
 e^{ \frac{x^3+1}{\sqrt{{pn}}} }
          % \binom{n}{m+j}
           %\pi(j)p^{j}q^{-j}\,.
%           \Big(\frac{p}{q}\Big)^{j}
\Bigg(
  \frac{1}{\sqrt{{ npq} }}
  +
  %\int_0^{ x\sqrt{{npq}}}
 %\sqrt{{1+\frac1n} }
  \int_0^{ x}
  e^{- \frac{s^2}{2} }
   %\exp\Big(- \frac{t^2}{2(n+1)pq} \Big)
   \,dt
  \frac{ds}{\sqrt{{ 2\pi} }}
   \Bigg)
           \,.
\end{equation}
\subsection{Lower bound.}
%Let us fix $j\geq 1$. To prepare for the lower bound, we rewrite 
Let us fix $j\in\zun$. To prepare for the lower bound, we rewrite 
$\pi(j)$ as follows:
\begin{multline}
  \label{zinop}
  \pi(j) \,=\, 
           \prod_{1\leq\ell\leq j}
  \frac{npq-(\delta+\ell-1)p}{npq+(\delta+\ell)q}
   \,\geq\, \prod_{1\leq\ell\leq j}
  \frac{npq-\ell p +pq}{npq+\ell q+pq}\cr
   \,\geq\, \prod_{1\leq\ell\leq j}
  \frac{(n+1)pq-\ell p }{(n+1)pq+\ell q}
   \,\geq\, \prod_{1\leq\ell\leq j}
   \frac{1}{1+
   \frac{\ell}{(n+1)pq-\ell p}}
  \,.
\end{multline}
We apply next the inequality~\eqref{lowi}.
For $j$ such that 
   $0\leq j\leq (n+1)pq/2$, we have
   \begin{gather}
   {(n+1)pq-j p}
   \,\geq\,
   (n+1)pq\Big(1-\frac{p}{2}\Big)
   \,\geq\,
   \frac12(n+1)pq\,\geq\,0\,,\cr
   \frac{j}{(n+1)pq-j p}\,\leq\,1\,,
   \end{gather}
  whence
\begin{comment}
  \label{zjnop}
  \pi(j) \,\geq \, 
  %         \prod_{1\leq\ell\leq j}
   %\exp\Big(-\frac{\ell}{(n+1)pq-\ell p}\Big)
   %\,=\, 
   \exp\Big(-
           \sum_{1\leq\ell\leq j}
   \frac{\ell}{(n+1)pq-\ell p}\Big)
   \cr
   \,\geq\,
   \exp\Big(-
           \sum_{1\leq\ell\leq j}
   \frac{\ell}{(n+1)pq-j p}\Big)
   \,=\,
   \exp\Big(-
   \frac{j(j+1)}{2\big((n+1)pq-j p\big)}\Big)
  \,.
\end{comment}
\begin{comment}
  \label{zjnop}
  \pi(j) \,\geq \, 
  %         \prod_{1\leq\ell\leq j}
   %\exp\Big(-\frac{\ell}{(n+1)pq-\ell p}\Big)
   %\,=\, 
   \exp\Big(-
           \sum_{1\leq\ell\leq j}
   \frac{\ell}{(n+1)pq-\ell p}\Big)
   \,\geq\,
   \exp\Big(-
   \frac{j(j+1)}{2\big((n+1)pq-j p\big)}\Big)
  \,.
\end{comment}
\begin{multline}
  \label{zjnop}
  \pi(j) \,\geq \, 
  %         \prod_{1\leq\ell\leq j}
   %\exp\Big(-\frac{\ell}{(n+1)pq-\ell p}\Big)
   %\,=\, 
   \exp\Big(-
           \sum_{1\leq\ell\leq j}
   \frac{\ell}{(n+1)pq-\ell p}\Big)\cr
   \,\geq\,
   \exp\Big(-
           \sum_{1\leq\ell\leq j}
   \frac{\ell}{(n+1)pq-j p}\Big)
   \,\geq\,
   \exp\Big(-
   \frac{j(j+1)}{2\big((n+1)pq-j p\big)}\Big)
  \,.
\end{multline}
Next, we use the inequality~\eqref{ine2}
\begin{comment}
  \forall t\in\Big[0,\frac{1}{2}\Big]\qquad 
  \frac{1}{1-t}\,\leq\, 1+2t
\end{comment}
to obtain that, for $0\leq j\leq (n+1)pq/2$,
\begin{multline}
  \label{nbnos}
  \frac{j(j+1)}{2\big((n+1)pq-jp\big)}
   \,=\,
\frac{j(j+1)}{2(n+1)pq}
   %\frac{1}{ \displaystyle 1+\frac{j}{(n+1)p}
   \frac{1}{ 1-\frac{j}{(n+1)q}
   }\cr
   \,\leq\,
\frac{j(j+1)}{2(n+1)pq}
   \Big( 1+\frac{2j}{(n+1)q}\Big)
   \,\leq\,
\frac{(j+1)^2}{2npq}
   +\frac{(j+1)^3}{n^2pq^2}
  \,.
\end{multline}
Substituting~\eqref{nbnos} in~\eqref{zjnop},
%and~\eqref{vanop}, we conclude that
we conclude that
\begin{equation}
  \label{ginos}
\forall j\,,\quad 0\leq j\leq \frac12(n+1)pq\,,\qquad
  \pi(j) 
  \,\geq \, 
   \exp\Big(-
\frac{(j+1)^2}{2npq}
-\frac{(j+1)^3}{n^2pq^2}
   \Big)
  \,.
\end{equation}
\begin{comment}
Contrary to the previous case, this inequality is as good as the one 
obtained by 
Feller 
(see \cite{F}, page 180, between formulas~$(2.3)$ and $(2.4)$).
\end{comment}
%\subsection{Upper bound of $P_{2n}(x)$.}
With the lower bound of $\pi(j)$ in hand, we will compute 
a lower bound on 
$P_{n}(x)$. Let us fix $x>0$ and let $n$ be an integer 
such that 
%$n\geq (q/p)x^2$. 
%$x\sqrt{{npq}}\geq 1$.
$x\sqrt{{npq}}\geq \max(1,2x^2)$.
This ensures that 
\begin{equation}
  \label{klfnst}
           {{0\leq j< x\sqrt{{npq}}}}
           \quad
           \Longrightarrow
           \quad
           %{j}\leq (n+1)p\,,\quad
j\leq 
\frac12(n+1)pq\,,\quad
         %  \frac{j}{n}\leq\frac12\,,\quad
   \frac{(j+1)^3}{n^2pq^2}
   %\leq\frac{x^3}{2}\sqrt{\frac{q}{pn}}\,.
   \leq\frac{8x^3}{\sqrt{{qn}}}\,.
\end{equation}
In particular, for these values of $j$, the inequality~\eqref{ginos}
holds, thus
\begin{equation}
  \label{srvinos}
           %{{0\leq j< x\sqrt{{npq}}}}
           %\quad
           %\Longrightarrow
           %\quad 
  \pi(j) 
  \,\geq \, 
   \exp\Big(-
\frac{(j+1)^2}{2npq}
   %+\frac{x^3}{2}\sqrt{\frac{q}{pn}}
   -\frac{8x^3}{\sqrt{{qn}}}
   \Big)
  \,\geq \, 
 e^{ -\frac{8x^3}{\sqrt{{qn}}} }
\int_{j+1}^{j+2} 
  e^{- \frac{t^2}{2npq} }
   %\exp\Big(- \frac{t^2}{2(n+1)pq} \Big)
   \,dt
  \,.
\end{equation}
Substituting inequality~\eqref{srvinos}
into the sum appearing in~\eqref{bwex}, we get
\begin{multline}
  \label{sawex}
  P_{n}(x)
           \,\geq\, 
           P(S_{n}=m)
          \,
 e^{ -\frac{8x^3}{\sqrt{{qn}}} }
  \int_1^{ x\sqrt{{npq}}+1}
  e^{- \frac{t^2}{2npq} }
   \,dt\cr
           \,\geq\, 
           %\sum_{\substack{k\in\zun\\ |k-n|< x\sqrt{\frac{n}2}}}
           %\frac{2}{2^{2n}} \binom{2n}{n}
           P(S_{n}=m)
   %\exp\Big( \frac{x^3}{2}\sqrt{\frac{q}{pn}}\Big)
 %e^{ \frac{x^3}{2}\sqrt{\frac{q}{pn}}}
           \,
 e^{ -\frac{8x^3}{\sqrt{{qn}}} }
          % \binom{n}{m+j}
           %\pi(j)p^{j}q^{-j}\,.
%           \Big(\frac{p}{q}\Big)^{j}
\Bigg(
  \int_0^{ x\sqrt{{npq}}}
  e^{- \frac{t^2}{2npq} }
   %\exp\Big(- \frac{t^2}{2(n+1)pq} \Big)
   \,dt
   -1
   \Bigg)
           \,.
\end{multline}
We make the change of variable 
$s { \sqrt{{npq}} }=t$ in~\eqref{sawex} and we get
\begin{equation}
  \label{aogex}
  P_{n}(x)
           \,\geq\, 
           %\sum_{\substack{k\in\zun\\ |k-n|< x\sqrt{\frac{n}2}}}
           %\frac{2}{2^{2n}} \binom{2n}{n}
           P(S_{n}=m)
   %\exp\Big( \frac{x^3}{2}\sqrt{\frac{q}{pn}}\Big)
 %e^{ \frac{x^3}{2}\sqrt{\frac{q}{pn}}}
 \,e^{ -\frac{8x^3}{\sqrt{{qn}}} }
          % \binom{n}{m+j}
           %\pi(j)p^{j}q^{-j}\,.
%           \Big(\frac{p}{q}\Big)^{j}
\Bigg(
  %\int_0^{ x\sqrt{{npq}}}
 \sqrt{{npq} }
  \int_0^{ x}
  e^{- \frac{s^2}{2} }
   %\exp\Big(- \frac{t^2}{2(n+1)pq} \Big)
   \,ds
   -1
   \Bigg)
           \,.
\end{equation}
We factorize 
$ \sqrt{{npq} }$,
we 
introduce the normalizing factor $\sqrt{2\pi}$, and we obtain
%and we use the inequality 
%%${ \sqrt{{1+\frac1n} }\leq e^{1/n}}$ to get
%$$\sqrt{{1+\frac1n} }\leq e^{\frac1n}
 %e^{ \frac{x^3+1}{\sqrt{{pn}}} }
%$$ 
%to get
\begin{equation}
  \label{bzoex}
  P_{n}(x)
           \,\geq\, 
           %\sum_{\substack{k\in\zun\\ |k-n|< x\sqrt{\frac{n}2}}}
           %\frac{2}{2^{2n}} \binom{2n}{n}
 \sqrt{{2\pi npq} }\,
           P(S_{n}=m)
   %\exp\Big( \frac{x^3}{2}\sqrt{\frac{q}{pn}}\Big)
 %e^{ \frac{x^3}{2}\sqrt{\frac{q}{pn}}}
 \,e^{- \frac{8x^3}{\sqrt{{qn}}} }
          % \binom{n}{m+j}
           %\pi(j)p^{j}q^{-j}\,.
%           \Big(\frac{p}{q}\Big)^{j}
\Bigg(
  %\int_0^{ x\sqrt{{npq}}}
 %\sqrt{{1+\frac1n} }
  \int_0^{ x}
  e^{- \frac{s^2}{2} }
   %\exp\Big(- \frac{t^2}{2(n+1)pq} \Big)
   \,dt
  \frac{ds}{\sqrt{{ 2\pi} }}
  -\frac{1}{\sqrt{{ npq} }}
   \Bigg)
           \,.
\end{equation}
The two non-asymptotic inequalities~\eqref{azoex} and~\eqref{bzoex} together imply the 
limit of de Moivre--Laplace in the case where $a=0$ and $b=x$. 
Applying the same result 
%with the variables $1-X_1,\dots,1-X_n$ and the parameter $1-p$, we 
with the parameter $1-p$, and setting 
\begin{equation}
\forall i\in\unn\qquad \tX_i=1-X_i\,,
\end{equation}
we obtain
\begin{equation}
  \label{tigd}
  \lim_{n\to\infty}\,P\Big(
    0<\frac{n-\big(\tX_1+\cdots+\tX_n\big)-n(1-p)}{\sqrt{np(1-p)}}<x
  \Big) \,=\,\frac{1}{\sqrt{2\pi}}
  \int_a^b
   e^{- \tfrac{s^2}{2} } ds\,.
  %\exp\big(-\frac{x^2}{2}\big)\,dx\,.
  %\int_a^b\exp\big(-\frac{x^2}{2}\big)\,dx\,.
\end{equation}
%Yet the variables $\tX_1,\dots,\tX_n$ are i.i.d. Bernoulli with parameter $p$, 
Yet $\tX_1,\dots,\tX_n$ are i.i.d. Bernoulli with parameter $p$.
Rewriting the probability as 
\begin{equation}
  \label{rigd}
  P\Big(
    0<\frac{-\big(\tX_1+\cdots+\tX_n\big)+np}{\sqrt{np(1-p)}}<x
  \Big) \,=\,
  P\Big(
    -x<\frac{\tX_1+\cdots+\tX_n-np}{\sqrt{np(1-p)}}<0
  \Big) 
  %\exp\big(-\frac{x^2}{2}\big)\,dx\,.
  %\int_a^b\exp\big(-\frac{x^2}{2}\big)\,dx\,.
\end{equation}
we see that
the limit~\eqref{tigd} corresponds to the case 
%obtain the limit for the case 
$a=-x$ and $b=0$. 
%Putting the two results together, we 
%obtain the limit for the case $a=-x$ and $b=0$. 
Putting the two previous results together, we obtain the limit 
for a symmetric interval.
%, and from there, as in the symmetric case, 
As in the symmetric case, this yields
the limit for an arbitrary interval $[a,b]$.
\bibliographystyle{vancouver}
\bibliography{clt.bib}

%\begin{biog}
%\item[Author Name 1] Insert author bio here.
%\begin{affil}
%Department of Mathematics, University America, Washington DC 20036\\
%authorname@ua.edu
%\end{affil}
%\end{biog}

%\begin{biog}
%\item[Author Name 2] Insert author bio here.
%\begin{affil}
%Department of Mathematics, University America, Washington DC 20036\\
%authorname@ua.edu
%\end{affil}
%\end{biog}
\begin{comment}
\item[Rapha\"el Cerf] is a former student of the Ecole Normale Supérieure of Paris.
He received his Ph.D. from the University of Montpellier in 1994. 
Since then, he has worked at 
%the Mathematics departement of
Universit\'e Paris--Sud,
except for the last ten years when he was seconded to the Ecole Normale Supérieure.
He is a probabilist and he spends desperately his mathematical energy trying to prove the famous 
conjecture 
$\theta(p_c,{\mathbb Z}^3)=0$ 
in percolation. He enjoys spending his time with his family, baking bread,
cooking, listening to classical music, biking and running.
\begin{affil}
LMO, Université Paris-Sud, CNRS, Université Paris--Saclay, 91405 Orsay.\\
Raphael.Cerf@universite-paris-saclay.fr
\end{affil}
\end{comment}

%\vfill\eject

\end{document}